\documentclass[arxiv, preprint, nameyear]{imsart}

\RequirePackage{amsthm,amsmath,amsfonts,amssymb}
\RequirePackage[colorlinks,citecolor=blue,urlcolor=blue]{hyperref}
\RequirePackage{graphicx}
\usepackage{multicol}

\usepackage{enumitem}
\usepackage{multirow}
\usepackage{graphicx}
\usepackage{array, makecell} 
\usepackage{diagbox}

\usepackage{natbib}

\startlocaldefs
\theoremstyle{plain}
\newtheorem{theorem}{Theorem}
\newtheorem{corollary}[theorem]{Corollary}

\newtheorem{remark}[theorem]{Remark}

\newcommand{\ep}{\textsf{E}}
\newcommand{\var}{\textsf{Var}}
\newcommand{\cov}{\textsf{Cov}}
\newcommand{\B}{\text{B}}
\renewcommand{\P}{\text{P}}
\newcommand{\iidsim}{\stackrel{\text{iid}}{\sim}}
\renewcommand{\d}{{\text{diff}}}
\renewcommand{\,}{{\text{, }}}

\newcommand{\varb}{\textsf{Var}_\text{B}}
\endlocaldefs
\begin{document}

\begin{frontmatter}
\title{Limiting distributions of graph-based test statistics on sparse and dense graphs}
%\title{A sample article title with some additional note\thanksref{t1}}
\runtitle{Limiting distributions of graph-based test statistics}
%\thankstext{T1}{A sample additional note to the title.}

\begin{aug}
%\author{\inits{F.}\fnms{Yejiong} \snm{Zhu}\ead[label=e1]{yjzhu@ucdavis.edu}} 
%\author{\inits{S.}\fnms{Hao} \snm{Chen}\ead[label=e2,mark]{hxchen@ucdavis.edu}}
\author{\fnms{Yejiong} \snm{Zhu}\ead[label=e1]{yjzhu@ucdavis.edu}}
\and
\author{\fnms{Hao} \snm{Chen}\ead[label=e2]{hxchen@ucdavis.edu}}
\runauthor{Y.~Zhu and H.~Chen}
%%%%%%%%%%%%%%%%%%%%%%%%%%%%%%%%%%%%%%%%%%%%%%
%% Addresses                                %%
%%%%%%%%%%%%%%%%%%%%%%%%%%%%%%%%%%%%%%%%%%%%%%
\address{University of California, Davis
%\printead{e1,e2}
}

\end{aug}

\begin{abstract}
Two-sample tests utilizing a similarity graph on observations are useful for high-dimensional and non-Euclidean data due to their flexibility and good performance under a wide range of alternatives. Existing works mainly focused on sparse graphs, such as graphs with the number of edges in the order of the number of observations, and their asymptotic results imposed strong conditions on the graph that can easily be violated by commonly constructed graphs they suggested. Moreover, the graph-based tests have better performance with denser graphs under many settings. In this work, we establish the theoretical ground for graph-based tests with graphs ranging from those recommended in current literature to much denser ones.
\end{abstract}


\end{frontmatter}

%%%%%%%%%%%%%%%%%%%%%%%%%%%%%%%%%%%%%%%%%%%%%%
%%%% Main text entry area:

\section{Introduction}

Given two random samples, $X_1,\cdots\, X_m \iidsim F_X,  Y_1,\cdots\, Y_n\iidsim F_Y$, we consider the hypothesis testing problem $H_0: F_X=F_Y$ against $ H_a: F_X\neq F_Y$. This two-sample testing problem is a fundamental problem in statistics and has been extensively studied for univariate and low-dimensional data. Nowadays, it is common that observations are in high dimensions \citep{cancer2012comprehensive,feigenson2014disorganization,zhang2020simple}, or non-Euclidean, such as networks \citep{bullmore2009complex,biswal2010toward,beckmann2021downregulation}. In many of these applications, one has little knowledge on $F_X$ or $F_Y$, making parametric tests unapproachable.

Nonparametric methods play important roles in solving two-sample testing problems and have a long history. For univariate data, some common choices are the Kolmogorov–Smirnov test \citep{smirnov1939estimation}, the Wald–Wolfowitz runs test \citep{wald1940test} and the Mann–Whitney rank-sum test \citep{mann1947test}. Since the middle of the 20th century, researchers have tried to extend these methods to multivariate data \citep{weiss1960two,bickel1969distribution}. The first practical test that can be applied to data in an arbitrary dimension or non-Euclidean data was proposed by \cite{friedman1979}, which is based on a similarity graph and is the start of graph-based tests. Over the years, graph-based tests evolved a lot and showed good power for a variety of alternatives and different kinds of data \citep{schilling1986multivariate,henze1999multivariate,rosenbaum2005exact,chen2013graph,generalized,chen2018weighted,chu2019asymptotic,zhang2020graph}. In the following, we give a brief review of the graph-based tests and discuss their limitations.

\subsection{A review of graph-based tests}\label{review}
\cite{friedman1979} proposed to pool all observations from both samples to construct the minimum spanning tree (MST), which is a tree connecting all observations such that the sum of edge lengths that are measured by the distance between two endpoints is minimized. They then count the number of edges that connect observations from different samples and reject $H_0$ when this count is significantly small. The rationale is that when two samples are from the same distribution, they are well mixed and this count shall be relatively large, so a small count suggests separation of the two samples and rejection of $H_0$. We refer this test to be the $\emph{\text{original edge-count test}}$ (OET). This test is not limited to the MST. \cite{friedman1979} also applied it to the $K$-MST\footnote{A $K$-MST is the union of the $1$st, $\cdots$ $K$th MSTs, where the 1st MST is the MST and the $k$th $(k> 1)$ MST is a tree connecting all observations that minimizes the sum of distance across edges subject to the constraint that this tree does not contain any edge in the 1st, $\cdots$ $k-1$th MST(s).}. Later, \cite{schilling1986multivariate} and \cite{henze1988multivariate} applied it to the $K$-nearest neighbor graphs ($K$-NNG), and \cite{rosenbaum2005exact} applied it to the cross-match graph.  \cite{zhang2020graph} extended the test to accommodate data with repeated observations.

Recently, \cite{generalized} noticed an issue of OET caused by the curse of dimensionality. They made use of a common pattern under moderate to high dimensions and proposed the $\emph{\text{generalized edge-count test}}$ (GET), which exhibits substantial power improvements over OET under a wide range of alternatives. Later, two more edge-count tests were proposed, the $\emph{\text{weighted edge-count}}$ $\emph{test}$ (WET) \citep{chen2018weighted} and the $\emph{max-type edge-count test}$ (MET) \citep{chu2019asymptotic}. WET addresses an issue of OET under unequal sample sizes, but it focuses on the locational alternatives. MET performs similarly to GET while it has some advantages under change-point settings.

In the following, we express the four graph-based test statistics with rigorous notations. The two samples $X_1, \cdots\, X_m$ and $Y_1, \cdots\, Y_n$, are pooled together and indexed by $1, \cdots\, N\ (N=m+n)$. Let $G$ be the set of all edges in the similarity graph, such as the $K$-MST. For an edge $e\in G$, let $e^+,e^-$ be two endpoints of the edge $e$, i.e. $e = (e^+,e^-)$.  Let $g_i$ be the group label of $i$-th observation with
\begin{eqnarray*}
    g_i=\left\{
\begin{aligned}
1 & \text{ if observation $i$ is from sample X} \\
2 & \text{ if observation $i$ is from sample Y,}
\end{aligned}
\right.
\end{eqnarray*}
and $R_1,R_2$ be the numbers of within-sample edges of sample X and sample Y, respectively, formally defined as
\begin{eqnarray*}
    R_j = \sum_{e\in G}1_{\{J_e=j\}},\quad j=1,2,
\end{eqnarray*}
where $1_{\{A\}}$ is the indicator function that takes value $1$ if event $A$ occurs and takes value $0$ otherwise, and
\begin{eqnarray*}
    J_e=\left\{
\begin{aligned}
0 & \quad\text{ if } g_{e^+}\neq g_{e^-} \\
1 & \quad\text{ if } g_{e^+}= g_{e^-} =1 \\
2 & \quad\text{ if } g_{e^+}= g_{e^-} =2.
\end{aligned}
\right.
\end{eqnarray*}
Since no distributional assumption was made for $F_X$ and $F_Y$, we use the permutation null distribution, which places probability $1 \slash { N\choose m}$ on each selection of $m$ observations among pooled observations as sample X. Let $\ep_\P, \var_\P, \cov_\P$ be the expectation, variance and covariance under the permutation null distribution.

{The four graph-based test statistics mentioned above can be expressed as follows:}
\begin{itemize}
    \item [1.] OET: $Z_o^\P=\big(R_1+R_2-\ep_\P(R_1+R_2)\big)/\sqrt{\var_\P(R_1+R_2)};$
    \item [2.] GET: $ S = \begin{pmatrix}R_1-\ep_\P(R_1)& R_2-\ep_\P(R_2) \end{pmatrix}\Sigma_R^{-1} \begin{pmatrix}R_1-\ep_\P(R_1) \\R_2-\ep_\P(R_2)\end{pmatrix},$ where $\Sigma_R = \var_\P(R_1,R_2)^T$;
    \item [3.] WET: $Z_w^\P = \big(R_w-\ep_\P(R_w))/\sqrt{\var_\P(R_w)},$ where $R_w = R_1(n-1)/(N-2)+R_2(m-1)/(N-2)$;
    \item [4.] MET: $\max\{Z_w^\P,Z_\d^\P\}$, where $Z_{\d}^\P = \big(R_{\d}-\ep_\P(R_\d)\big)/\sqrt{\var_\P(R_\d)}$ with $R_\d= R_1-R_2$.
\end{itemize}
The analytic formulas for the expectations and variances of $R_1+R_2$, $(R_1,R_2)^T$, $R_w$, $R_\d$ can be found in  \cite{generalized, chen2018weighted, chu2019asymptotic}. It was also shown in \cite{chu2019asymptotic} that the statistic $S$ can be decomposed as
\begin{equation}
    S = (Z_w^\P)^2+(Z_{\d}^\P)^2 \label{decom_S}\\
\end{equation}
and
$\cov_\P(Z_w^\P,Z_\d^\P)=0.$

\subsection{Limitations of existing theorems on graph-based test statistics}\label{restriction_on_cur_theorem}
Since the computation of $p$-values by drawing random permutations is time-consuming, approximated $p$-values from the asymptotic null distributions of the graph-based test statistics are useful in practice. Some existing theorems provide sufficient conditions for the validity of the asymptotic distributions. Before stating these results, we first define some essential notations. 

\begin{figure}
    \centering
    \includegraphics[width = 0.35\textwidth]{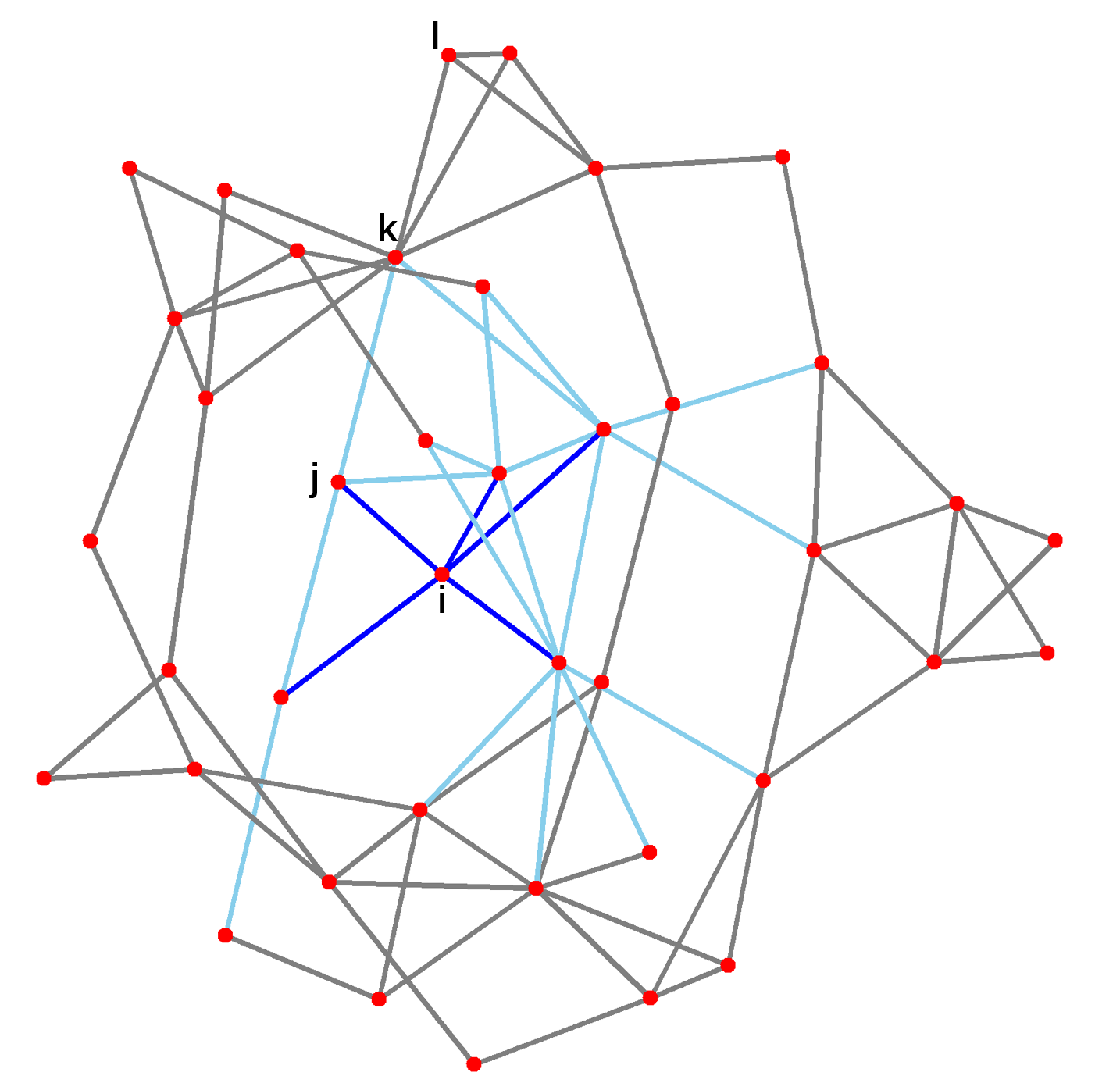}
    \includegraphics[width = 0.35\textwidth]{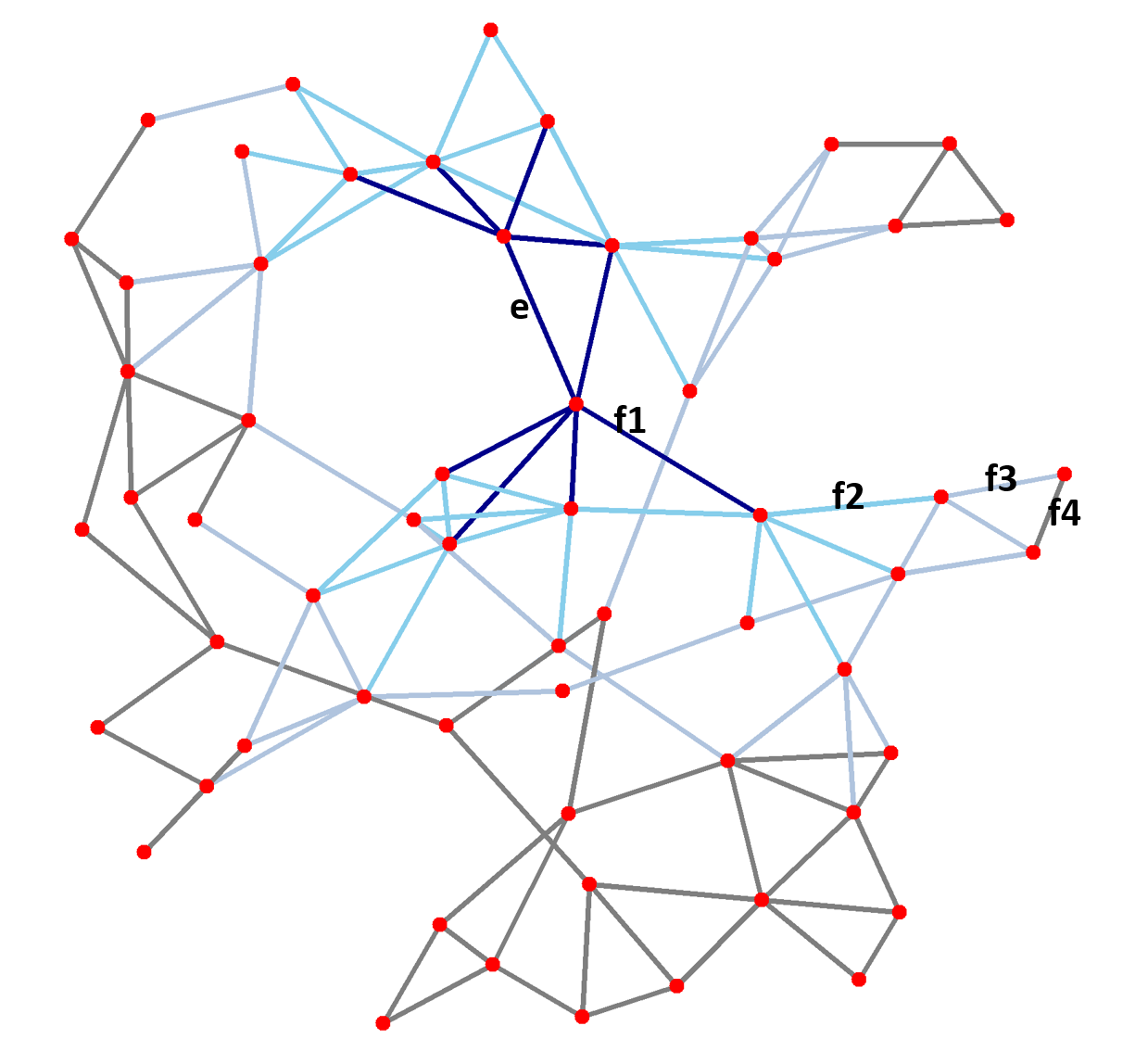}
    \caption{The left panel is an illustration of notations related to a node $i$: $G_i=\{\text{blue edges}\}$, $G_{i,2} = \{\text{sky blue edges}\}\cup\{\text{blue edges}\}$; $j\in node_{G_i},k\in node_{G_{i,2}}\backslash node_{G_i},l\notin node_{G_{i,2}}$, and the right panel is an illustration of notations related to an edge $e$: $ A_e= \{\text{dark blue edges}\}, B_e\backslash A_e=$ \{\text{sky blue edges}\}, $C_e\backslash B_e =$ \{light steel blue edges\}; $f_1\in A_e, f_2\in B_e\backslash A_e, f_3\in C_e\backslash B_e, f_4\notin C_e $. Both graphs are plotted by the `ggnet2' function in R.}
    \label{notation_graph}
\end{figure}

For each node $i\in \mathcal{N} \triangleq \{1, \cdots\, N\}$, let $G_i$ be the set of edges with one endpoint node $i$, $node_{G_i}$ be the set of nodes connected by $G_i$ excluding node $i$, $G_{i,2}$ be the set of edges with at least one endpoint in $node_{G_i}$, and $node_{G_{i,2}}$ be the set of nodes connected by $G_{i,2}$ excluding node $i$. For each edge $e=(e^+,e^-)\in G$, define $A_e = G_{e^+}\cup G_{e^-}$, $B_e = G_{e^+,2}\cup G_{e^-,2}$ and $C_e$ be the set of edges that share at least one common node with an edge in $B_e$. We use $|\cdot|$ to denote the cardinarity of a set. Then $|G_i|$ is the degree of the node $i$. Figure \ref{notation_graph} plots the quantities related to a node $i$ and an edge $e$, respectively.

We also define $\Tilde{d}_i = |G_i|-2|G|/N$ to be the centered degree of node $i$ and $V_G = \sum_{i=1}^N\Tilde{d}_i^2=\sum_{i=1}^N|G_i|^2-4|G|^2/N$ that measures the variability of $|G_i|$'s.  Besides, $a_n= o(b_n)$ or $a_n\prec b_n$ means that $a_n$ is dominated by $b_n$ asymptotically, i.e. $\lim_{n\rightarrow \infty}a_n/b_n=0$, $a_n\precsim b_n$ means that $a_n$ is bounded above by $b_n$ (up to a constant factor) asymptotically, and $a_n= O(b_n)$ or $a_n \asymp b_n$ means that $a_n$ is bounded both above and below by $b_n$ asymptotically. We use $a\wedge b$ for $\min\{a,b\}$. For two sets $S_1$ and $S_2$, $S_1\backslash S_2$ is used for the set that contains elements in $S_1$ but not in $S_2$.

\begin{table}\label{exiting th}
\caption{ Major existing works on graph-based tests and their conditions on the graph for  asymptotic distributions.  }
\label{table1}
\centering
\begin{tabular}{llll}
\hline
 &Test Statistic& Graph conditions & \makecell[l]{max size \\of possible graphs} \\
 \hline
 \cite{friedman1979}&Original &MST with $\sum_{i=1}^N|G_i|^2 = O(N)$& $N-1$  \\
  \hline
\cite{schilling1986multivariate}&Original&\makecell[l]{$K$-NNG, $K=O(1)$\\ for low-dimensional data}& $O(N)$ \\
 \hline
 \cite{henze1988multivariate} &Original&\makecell[l]{$K$-NNG, $K=O(1)$\\ with bounded maximal indegrees} &$O(N)$
 \\\hline
 \cite{rosenbaum2005exact}&Original &cross-match graph & $N/2$ or $(N-1)/2$ \\
 \hline
\cite{chen2015graph} &Original &
    \makecell[l]{$|G| = O(N^{\alpha})\, 0<\alpha<1.125$\\ $\sum_{e\in G}|A_e||B_e| = o(N^{1.5(\alpha\wedge 1)})$}
    
& $O(N^\alpha)\,\alpha<1.125$ \\
\hline
\cite{generalized}&Generalized& \makecell[l]{$|G| = O(N)$\\$\sum_{i=1}^N|G_i|^2 = O(N)$\\ $\sum_{e\in G}|A_e||B_e| = o(N^{1.5})$}& $O(N)$   \\
 \hline
\cite{chen2018weighted} &Weighted & \makecell[l]{$|G| = O(N^\alpha)\, 1\leq\alpha<1.25$\\$\sum_{e\in G}|A_e||B_e| = o(N^{1.5\alpha})$\\$\sum_{e\in G}|A_e|^2 = o(N^{\alpha+0.5})$} & $O(N^\alpha)$, $\alpha < 1.25$\\
\hline
\cite{chu2019asymptotic}& \makecell[l]{Generalized \\ Weighted \\Max-type} & \makecell[l]{$|G| = O(N^\alpha), 1\leq\alpha<1.25$\\$\sum_{e\in G}|A_e||B_e| = o(N^{1.5\alpha})$\\$\sum_{e\in G}|A_e|^2 = o(N^{\alpha+0.5})$\\$V_G = O(\sum_{i=1}^N|G_i|^2)$}&  $O(N^\alpha)\,  \alpha < 1.25$\\
\hline
\multicolumn{4}{l}{\footnotesize{\makecell[l]{$^*$For the conditions in \cite{chen2018weighted} and \cite{chu2019asymptotic}, the size of the graph is bounded\\ by the condition on $|A_e|^2$:
$\sum_{e\in G}|A_e|^2 = o(N^{\alpha+0.5})$ requires that $\alpha<1.25$.}}}
\end{tabular}
\end{table}

\begin{figure}
\includegraphics[width=0.38\textwidth,height = 0.3\textwidth]{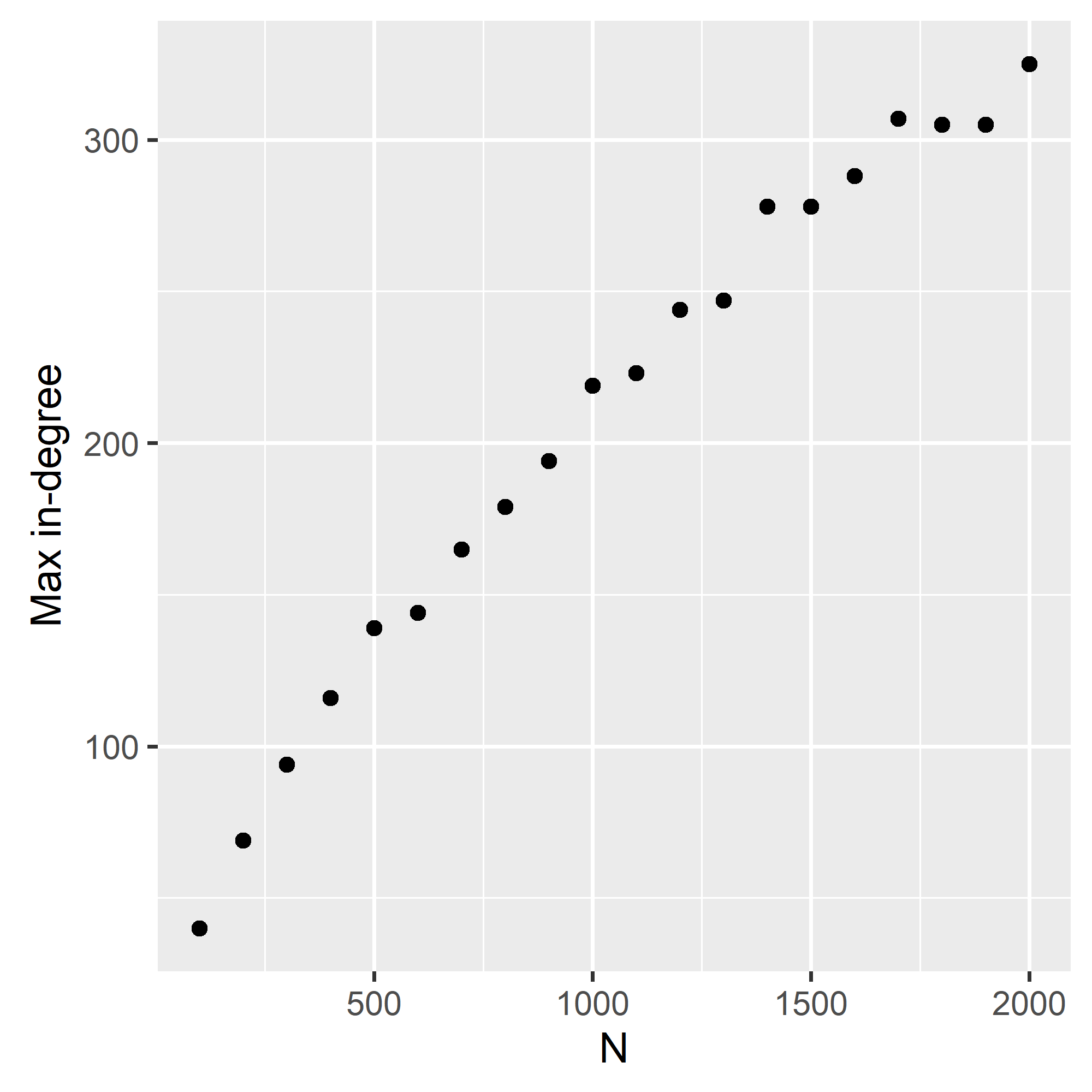}
\includegraphics[width=0.38\textwidth,height = 0.3\textwidth]{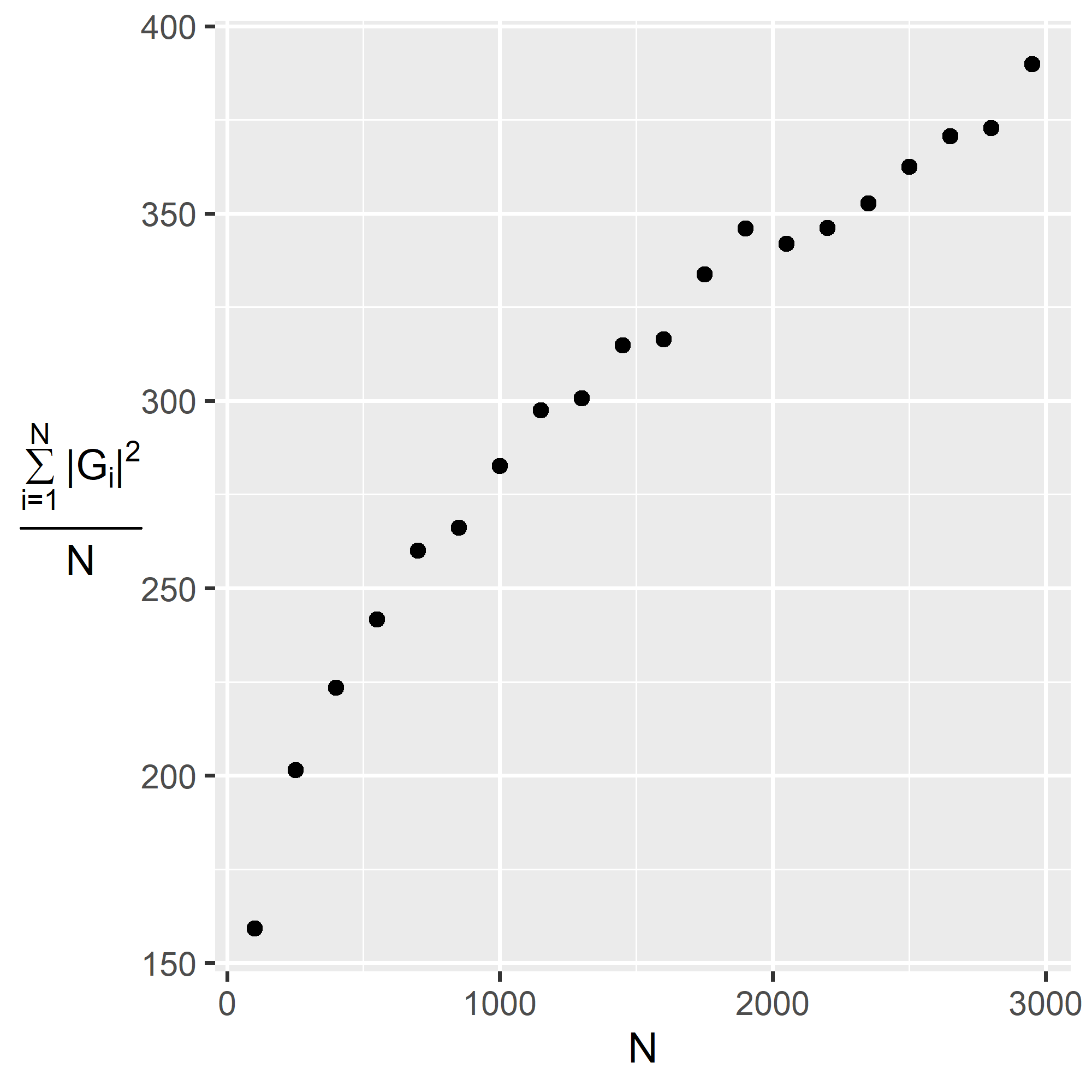}
\includegraphics[width=0.38\textwidth,height = 0.3\textwidth]{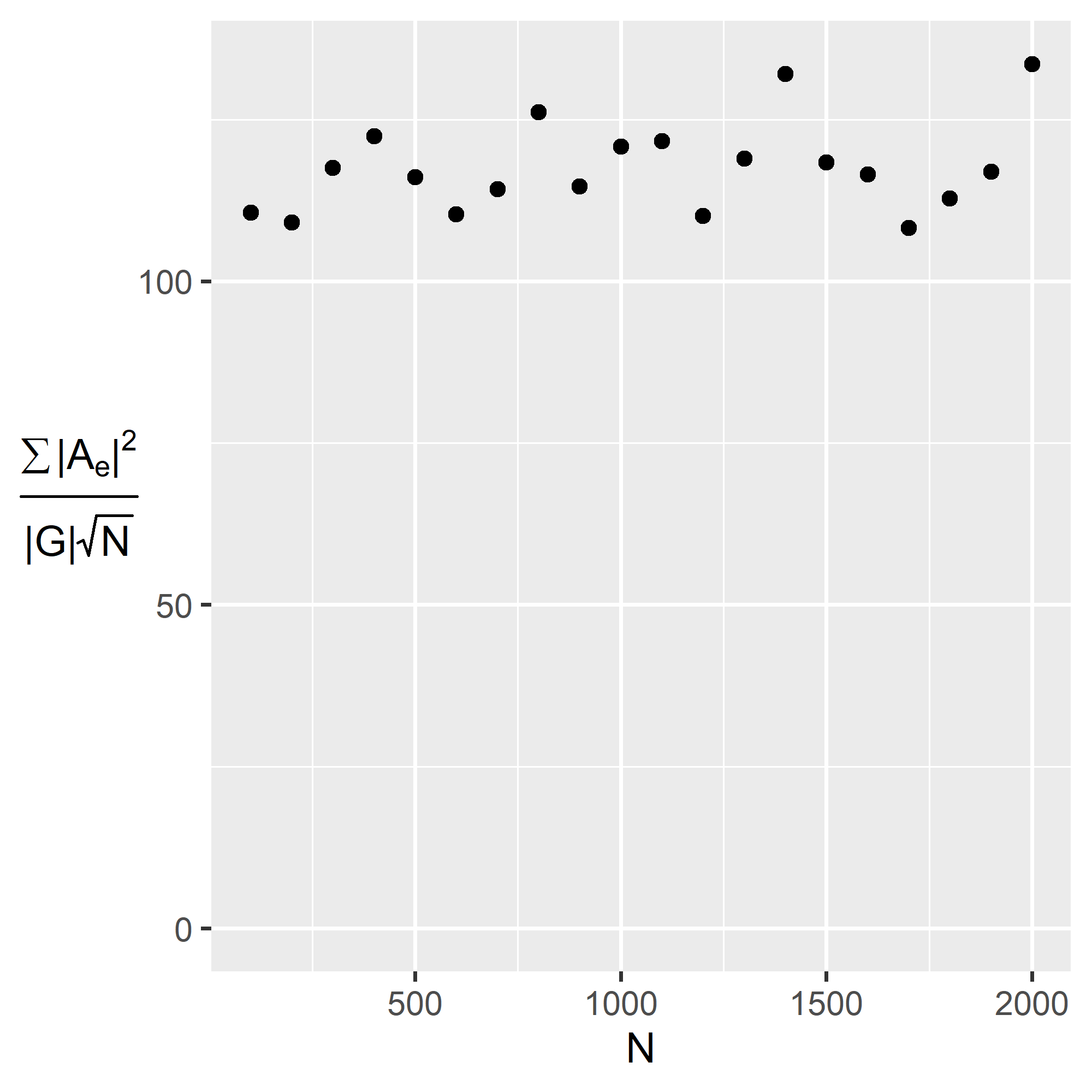}
\includegraphics[width=0.38\textwidth,height = 0.3\textwidth]{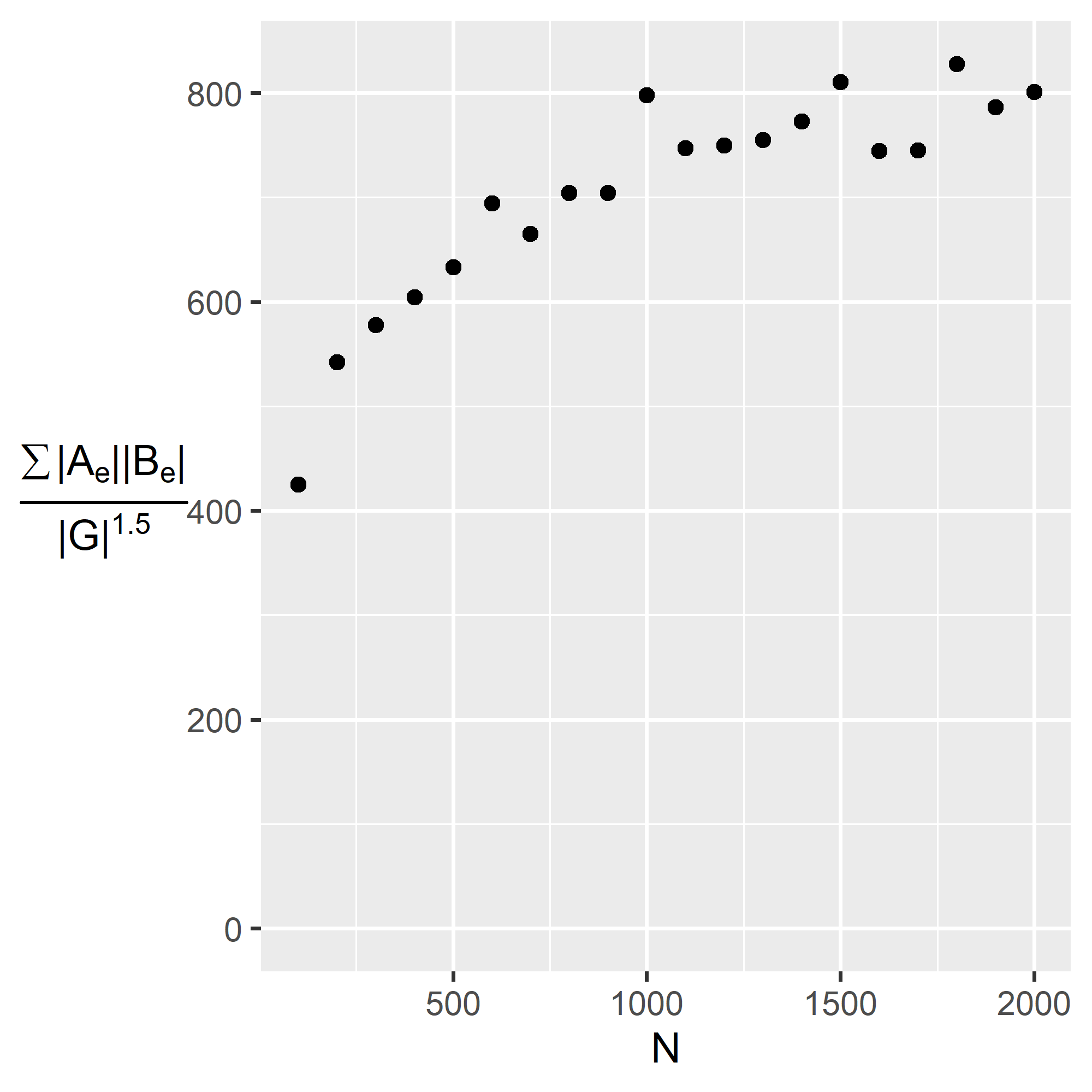}
\caption{\label{condition_existing} Key quantities with respect to $N$. Top-left: the maximum in-degree in the $5$-NNG;  top-right: the ratio of $\sum_{i=1}^N|G_i|^2$ to $N$ in the 5-MST; bottom-left: the ratio of $\sum |A_e|^2$ to $|G|\sqrt{N}$ in the $5$-MST; bottom-right: the ratio of $\sum |A_e||B_e|$ to $|G|^{1.5}$ in the 5-MST.}

\end{figure}

The major existing works that studied the asymptotic null distributions of the graph-based test statistics are listed in Table \ref{table1}. In general, these works put requirements on the graph such as the maximum in-degree, $\sum |G_i|^2$, $\sum |A_e|^2$ and $\sum |A_e||B_e|$. The conditions in \cite{friedman1979} are limited to the MST, while $K$-MST with $K>1$ has a better performance in general \citep{generalized}. For those that are more relaxed on the graph and data, \cite{henze1988multivariate} requires a bounded maximal in-degree in $K$-NNG, \cite{generalized} requires $\sum |G_i|^2$ to be of the same order as $N$, and \cite{chen2018weighted} and \cite{chu2019asymptotic} provide the current weakest conditions that require $\sum |A_e|^2 = o(|G|\sqrt{N})$ and $\sum |A_e||B_e| = o(|G|^{1.5})$. However, those conditions are often too strong to hold under even simple scenarios. For example, we generate two samples with equal sample size ($m=n=N/2$) from the 500-dimensional standard multivariate normal distribution, and construct the $5$-MST and $5$-NNG using the Euclidean distance. Figure \ref{condition_existing} plots the maximum in-degree, $\sum_{i=1}^N|G_i|^2/N$, $\sum |A_e|^2/(|G|\sqrt{N})$ and $\sum |A_e||B_e|/|G|^{1.5}$ with different $N$'s. We see that the conditions on them are badly violated: the maximum in-degree goes up rather than bounded by a constant, $\sum_{i=1}^N |G_i|^2/N$ increases with $N$ rather than bounded by a constant, and $\sum_{e\in G} |A_e|^2/(|G|\sqrt{N})$ and $\sum_{e\in G} |A_e||B_e|/|G|^{1.5}$ stay at a large value (in hundreds) as $N$ increases rather than $o(1)$ as required by the conditions.  If we make the graph denser, such as 10-MST, these conditions are even more badly violated. However, in many settings, the graph-based tests work better under denser graphs (see Section \ref{merit_denser}).

\subsection{The merits of denser graphs in improving power for graph-based tests}\label{merit_denser}
\cite{friedman1979} found that the original edge-count test in general had a higher power under the 3-MST than that under the 1-MST. Similarly, \cite{generalized} found that the generalized edge-count test in general had a higher power under the 5-MST than that under the 1-MST. We here check the performance of these tests under even denser graphs. In particular, for $m=n=100$, we consider the generalized edge-count tests on the $5$-MST $(\text{GET}_5)$ and on the $50$-MST $(\text{GET}_{50})$, the original edge-count tests on the $5$-MST $(\text{OET}_5)$ and on the $50$-MST $(\text{OET}_{50})$. All $K$-MSTs here are constructed under the Euclidean distance. We also include two other tests as baselines: the kernel two-sample test in \cite{gretton2012kernel} with the $p$-value approximated by 10,000 bootstrap samples (Kernel) and the Adaptable Regularized Hotelling's $\text{T}^2$ test \citep{li2020adaptable} (ARHT).

\begin{figure}
\includegraphics[width = 0.4\textwidth]{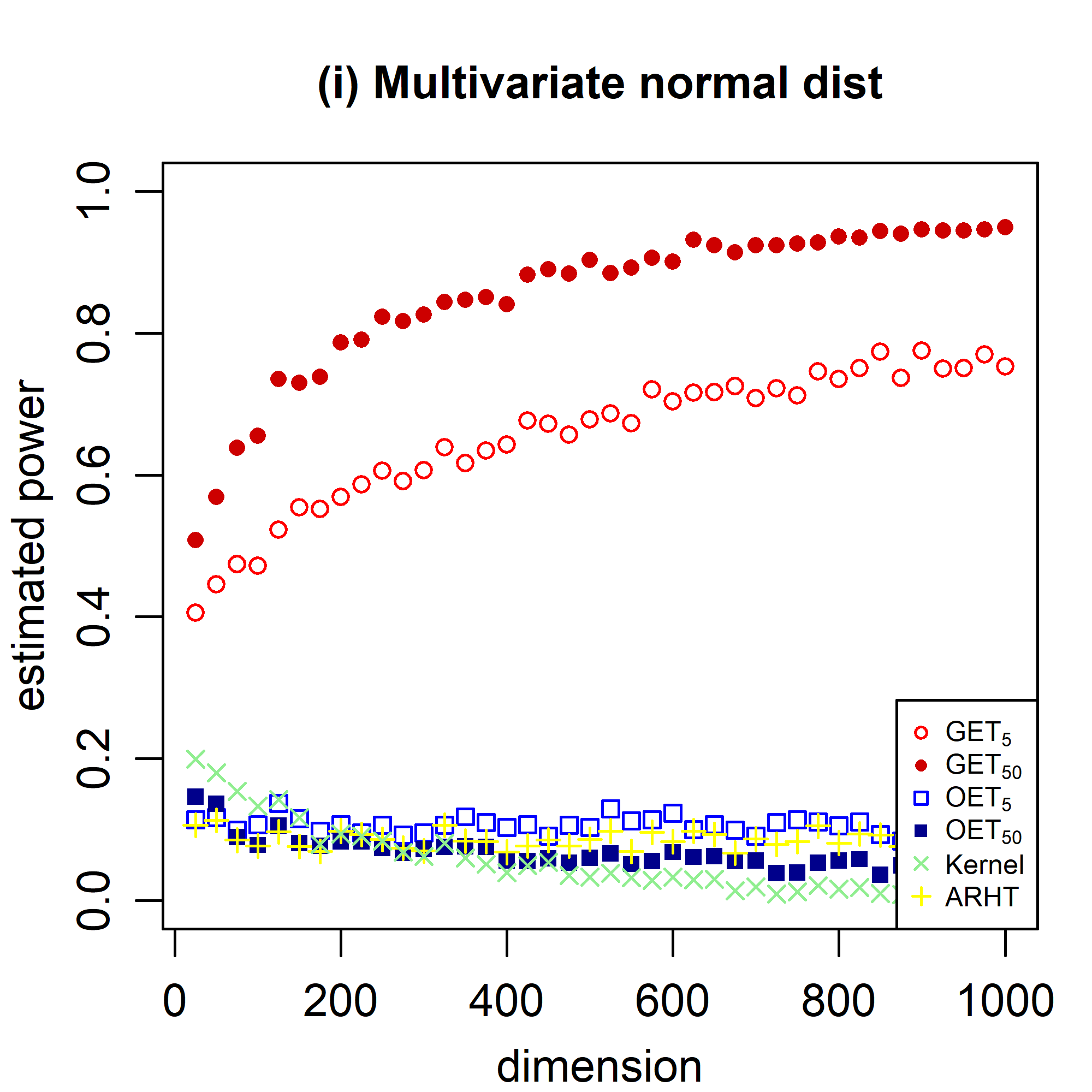}
\includegraphics[width = 0.4 \textwidth]{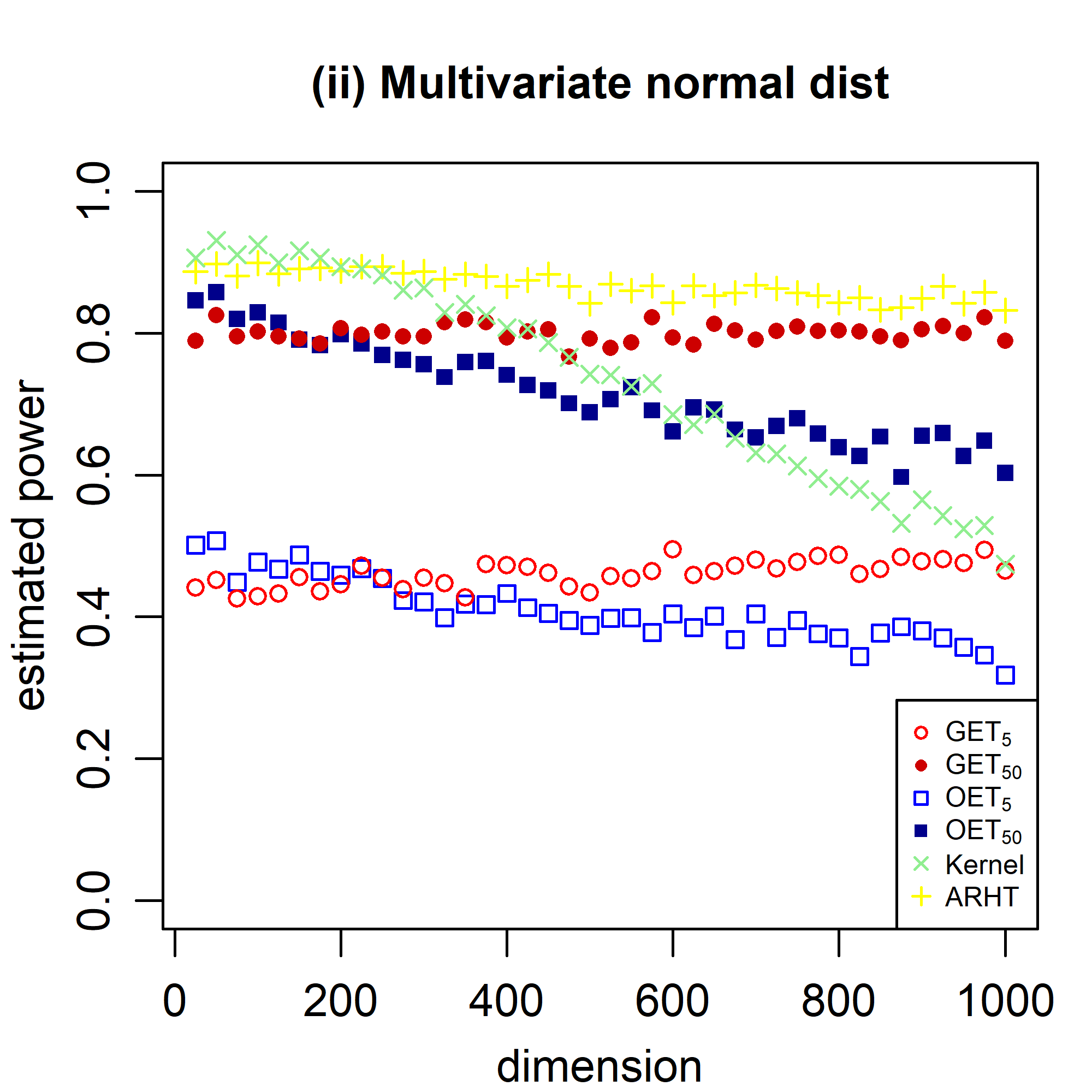}
\includegraphics[width = 0.4\textwidth]{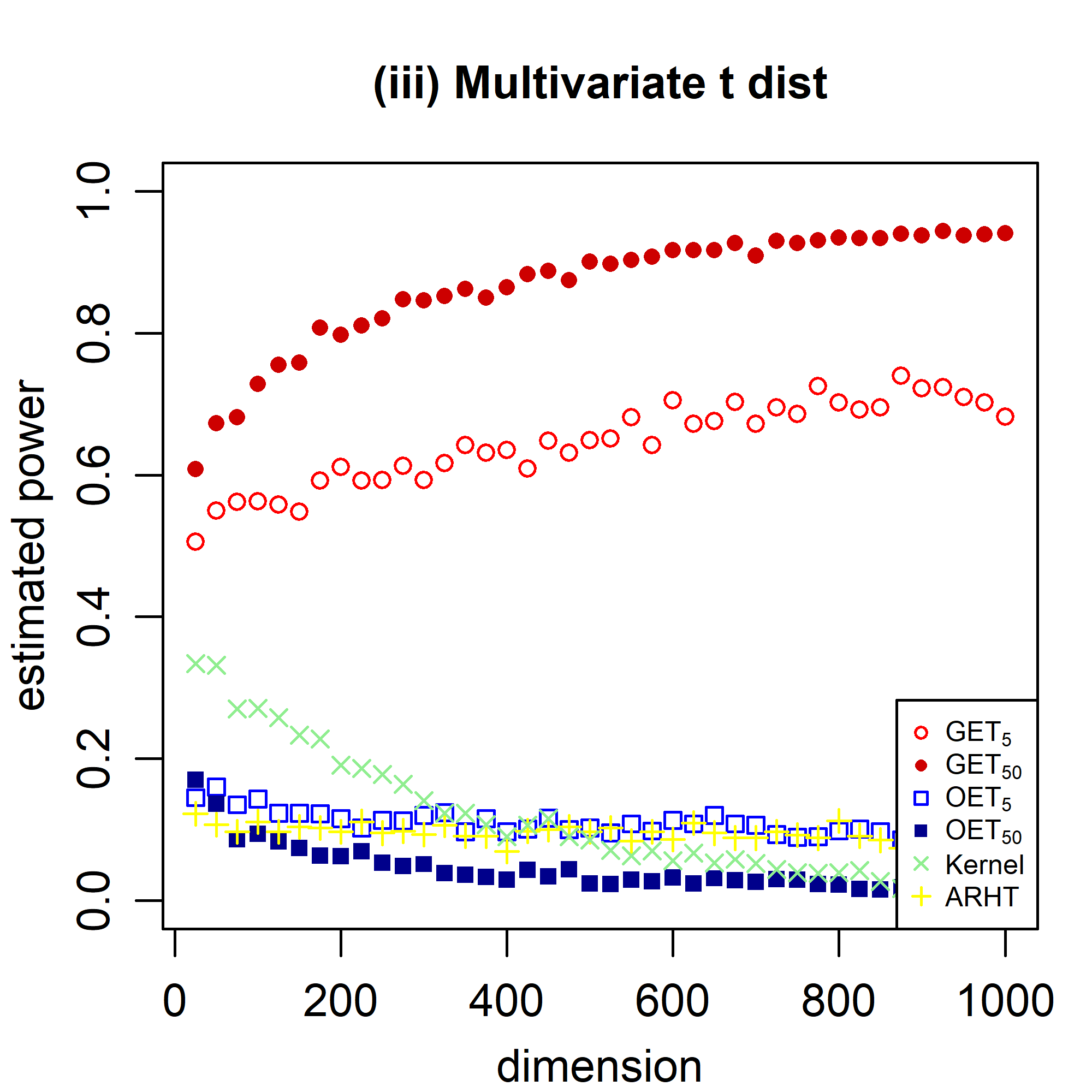}
\includegraphics[width = 0.4\textwidth]{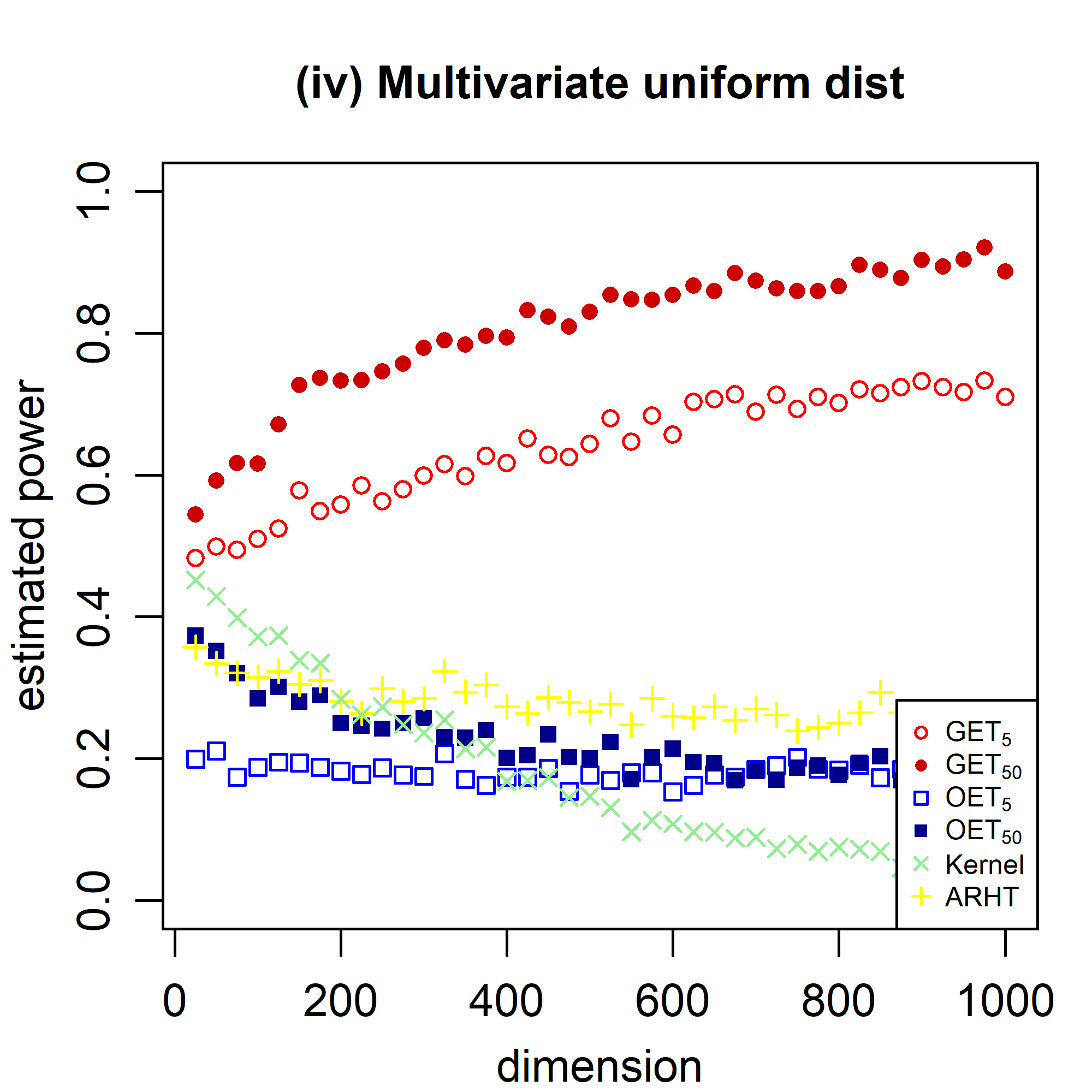}
\caption{Estimated power of the generalized edge-count tests on the $5$-MST $(\text{GET}_5)$ and on the $50$-MST $(\text{GET}_{50})$, the original edge-count tests on the $5$-MST $(\text{OET}_5)$ and on the $50$-MST $(\text{OET}_{50})$, the kernel test (Kernel), and the Adaptable Regularized Hotelling's $\text{T}^2$ test (ARHT) under different simulation settings (i) - (iv).}
\label{fig:Estimated_power}
\end{figure}

We consider different distributions in the comparison. Explicitly, 
\begin{align*}
    &X_i = \Sigma^\frac{1}{2}U_i,\text{ } i=1, \cdots \, 100,\quad Y_j = (1+ad^{-\frac{1}{3}})\Sigma^\frac{1}{2}V_j-bd^{-\frac{1}{3}}\mathbf{1}_d,\text{ }j=1, \cdots \, 100,
\end{align*}
with $\Sigma = (0.5^{|i-j|})_{1\leq i,\space j\leq d}$, where $d$ is the dimension of the data, $\mathbf{1}_d$ is a $d$-length vector of all ones. Let $\mathbf{0}_d$ be a $d$-length vector of all zeros and $I_d$ be a $d$-dimensional identity matrix. We consider four different settings:
\begin{itemize}
    \item[(i) ] $U_1,\cdots\,U_m,V_1,\cdots\,V_n \iidsim N(\mathbf{0}_d,I_d)$, $a=b=0.17$,
    \item[(ii) ] $U_1,\cdots\,U_m,V_1,\cdots\,V_n \iidsim N(\mathbf{0}_d,I_d)$, $a= 0.1, b=0.6$,
    \item[(iii) ] $U_1,\cdots\,U_m,V_1,\cdots\,V_n \iidsim t_5(\mathbf{0}_d,I_d)$, $a=b=0.25$,
    \item[(iv) ] $U_1,\cdots\,U_m,V_1,\cdots\,V_n \iidsim \text{Uniform} [-0.5,0.5]^d$, $a=0.12,b=0.1.$
    
\end{itemize}

Here, $a$ and $b$ are chosen so that the tests have moderate power in low dimensions. The dimension $d$ ranges from 25 to 1000 with an increment of 25. The power of tests are estimated through 1,000 simulation runs (Figure \ref{fig:Estimated_power}). We see that $\text{GET}_{50}$ works well in these settings, either having the best power or on par with the test of the best power.  $\text{GET}_5$ in general has a lower power than $\text{GET}_{50}$. For OET, it is only powerful under setting (ii). The worse performance of OET compared to GET is expected as OET covers less alternatives than GET for high-dimensional data \citep{generalized}. Under setting (ii) where OET is powerful, $\text{OET}_{50}$ has a higher power than $\text{OET}_{5}$.

%\begin{figure}
%    \centering
%    \includegraphics[scale=0.5]{powerwithk.png}
%    \caption{Estimated power of GET under the $K$-MST with different $K$ values.}
%    \label{fig:mixed}
%\end{figure}

\subsection{Our contribution}
From Section \ref{merit_denser}, we see that the use of denser graphs has a promising effect in improving power for graph-based tests. So far, the best theoretical results on dense graphs are in \cite{chen2018weighted} and \cite{chu2019asymptotic}, which allow the maximum size of possible graphs to be of order $N^\alpha$, $1\leq\alpha<1.25$. But in the numerical studies in Section \ref{merit_denser} where $N = 200$, the number of edges in the $50$-MST is $9950$ $(\approx N^{1.74})$. Existing conditions cannot work for such dense graphs. In addition, even for sparse graphs, current existing conditions usually do not hold (see Section \ref{restriction_on_cur_theorem}). Therefore, it is important to figure out whether the conditions on graphs can be weakened and to what extent. Throughout the paper, we consider simple undirected graphs that contain no duplicate edges and no loops.
% the relaxed conditions for the validity of the limiting distributions for both sparse graphs and denser graphs.

%\cite{friedman1979} used the moment-based method proposed in \cite{daniels1944relation}. Since then, the moment-based method has been extended. \cite{friedman1983graph} claimed that Daniels' conditions for asymptotic normality can be weakened and gave their relaxed conditions. However, they didn't give an explicit proof and merely stated that they could be obtained by a similar argument as in \cite{daniels1944relation}. \cite{pham1989asymptotic} found that conditions in \cite{friedman1983graph} are not sufficient. They fixed this problem and proposed a new set of relaxed conditions. By using these conditions, the results for the graph-based tests can be relaxed. However, we found that, our `locSCB' approach  could result in much more relaxed conditions for $S$. In the main context in this section, we focus on the `locSCB' approach to prove the Theorem \ref{th1}. The proof of Theorems \ref{th2}, \ref{th3} and \ref{th4} are similar to the proof of Theorem \ref{th1}, and are deferred to Supplement S1. Some comparisons to the moment-based method are provided in the Supplement S2.

\cite{friedman1979} applied the moment-based method in \cite{daniels1944relation} to derive sufficient conditions for the asymptotic normality of $Z_o^\P$, and \cite{generalized,chen2018weighted,chu2019asymptotic} made use of the bootstrap null distribution and the second neighbor dependent Stein's method to show the asymptotic normality of $(R_1,R_2)^T$ and $Z_w^\P$ under the permutation null distribution. In this paper, we seek improvements in both directions, especially the latter one. In particular, we propose to use a ``locSCB'' (local Stein's method on Conditioning Bootstrap) approach that use Stein's method to carefully deal with all the first neighbor dependency under the bootstrap null distribution and link the permutation null distribution and bootstrap null distribution through conditioning.  By doing this, we are able to weaken the conditions on the graphs to a tremendous amount. Under new conditions, the maximum size of possible graphs can be as large as $M^{1-\epsilon}$ with an arbitrarily small $\epsilon$, where $M = {N\choose 2}$ is the size of the complete graph. In addition, we also quantify the upper bound of Stein's inequality.

The ``locSCB'' approach is not limited to show the asymptotic properties of these four graph-based test statistics listed in Section \ref{review}. It can be applied to some other nonparametric two-sample test statistics and $K$-sample test statistics under the permutation null distribution, and to the change-point analysis settings.
The main theorems are provided in Section \ref{sec2}. We discuss the new conditions in Section \ref{some discussion on theorems} from various aspects. Section \ref{proof of th1} provides detailed proof for the theorem on GET. The detailed proofs for other edge-count tests are deferred to the Supplementary Material \citep{zhu2023}.

\section{Asymptotic distribution under denser graphs}\label{sec2}

 We first state the main results for the four statistics. The discussion of new conditions are deferred to Section \ref{some discussion on theorems}.  The proofs are provided in Section \ref{proof of th1} and the Supplementary Material \citep{zhu2023}.
 
 We use $\xrightarrow{\mathcal{D}}$ to denote convergence in distribution, and use `the usual limit regime' to refer $N\rightarrow \infty$ and $\lim_{N\rightarrow \infty} m/N = r\in (0,1)$. We define $N_{sq}$ as the number of squares in the graph, $N_{i,j}$ as the number of nodes connecting to nodes $i$ and $j$ simultaneously, and recall that $\Tilde{d_i} = |G_i|-2|G|/N$ (defined in Section \ref{restriction_on_cur_theorem}). { The conditions needed for the asymptotic results of the four statistics are listed below.}

\begin{enumerate}[label=\textbf{C.\arabic*}]
\item $\sum_{i=1}^N|G_i|^2 = o(|G|^{1.5})$, $N_{sq} = o(|G|^2)$; \label{cond for Zw}
\item $\sum_{i=1}^N \left|\Tilde{d_i}\right|^3 = o(V_G^{1.5})$, $\sum_{i=1}^N \Tilde{d_i}^3 = o(V_G\sqrt{|G|})$, $\sum_{i=1}^N\sum_{j,k\in node_{G_i}}^{j\neq k}\Tilde{d_j}\Tilde{d_k} = o(|G|V_G)$; \label{add cond for S}
\item $\max(\Tilde{d_i}^2) = o(V_G)$; \label{cond for Zd}
\item \label{cond for Zo} $\sum_{i=1}^N|G_i|^2 = o(T^{1.5})$, $\sum_{i=1}^N\left|\Tilde{d_i}\right|^3 = o(T^{1.5})$, $\sum_{i=1}^N\sum_{j,k\in node_{G_i}}^{j\neq k}\Tilde{d_j}\Tilde{d_k} = o(T^2)$, with $T = |G| + V_G$.
\end{enumerate}

\subsection{Generalized edge-count test}

\begin{theorem}\label{th1}
Under {Condition \ref{cond for Zw} and \ref{add cond for S}}, in the usual limit regime, 
$S\xrightarrow{\mathcal{D}} \chi_2^2$ under the permutation null distribution.
\end{theorem}

A detailed comparison of the conditions in Theorem \ref{th1} to the best existing conditions is provided in Section \ref{comparion of condition}. We will see that Theorem \ref{th1} provides much weaker conditions.

The conditions in Theorem \ref{th1} are not easily understandable for a graph. In the following, we provide a set of conditions that only involve up to the second moment of the degree distribution, and are easier to understand. Let $Q_G$ be a random variable generated from the degree distribution of a graph $G$ built on $N$ nodes. Then it is not hard to see that $\ep (Q_G) = 2|G|/N$ and $\var(Q_G) = V_G/N$.

\begin{corollary}\label{corollary 1}
Suppose $|G| =O(KN)$ with $1\precsim K \prec N$, if $\max(1,K^2/N)\precsim \var(Q_G) \prec K^{1.5}\sqrt{N}$, and the concentration inequality
\begin{align}
    P(|Q_G-\ep(Q_G)|\geq t)\leq 2\mathrm{exp}\left(-\frac{ct^2}{N^a}\right), \quad t>0, \label{concentration ineq}
\end{align}
holds for all large $N$ with some constants $c>0$ and $0<a<1$,
then,
in the usual limit regime, 
$   S\xrightarrow{\mathcal{D}} \chi_2^2$ under the permutation null distribution.

\end{corollary}

\begin{remark}\label{remark KMST for S}
\textnormal{Corollary \ref{corollary 1} is derived from Theorem \ref{th1} at the cost of sacrificing some of its generality. The conditions in Corollary \ref{corollary 1} could be further weakened when more information on the degree distribution is available. For the $K$-MST constructed on multivariate data, the maximum degree of the MST has been studied for fixed dimensions \citep{robins1994maximum}, which is of $O(1)$. Hence, under fixed dimensions, the maximum degree of the $K$-MST has the order at most of $K$, which is sufficient for the conditions in Corollary \ref{corollary 1} to hold when $K = O(N^\beta)$, $\beta<0.5$. {Similarly, under fixed dimensions, the maximum degree of the $K$-NNG also has the order of at most $K$, and thus the conditions in Corollary \ref{corollary 1} also hold for the $K$-NNG when $K = O(N^\beta)$, $\beta<0.5$.} To further relax $K$ and/or dimensions, the degree distributions of the $K$-MST {or the $K$-NNG} is needed, which is nontrivial for both fixed and non-fixed dimensions and will be explored in future research.}
\end{remark}

\subsection{Weighted and max-type edge-count test}

\begin{theorem}\label{th2} 
Under {Condition \ref{cond for Zw}}, in the usual limiting regime, 
$ Z_w^\P\xrightarrow{\mathcal{D}} N(0,1)$ under the permutation null distribution.
\end{theorem}

\begin{theorem}\label{th3}
Under {Condition \ref{cond for Zd}}, in the usual limit regime, 
$  Z_\d^\P\xrightarrow{\mathcal{D}} N(0,1) $ under the permutation null distribution.
\end{theorem}

Condition \ref{cond for Zd} is equivalent to
\begin{align} \label{equivalent condition to max type}
    \frac{\sum_{i=1}^N\left|\Tilde{d_i}\right|^{2+\delta}}{V_G^{\frac{2+\delta}{2}}} \rightarrow 0,\quad  \text{ for some } \delta >0,
\end{align} 
according to \cite{hoeffding1951combinatorial}. One condition in \ref{add cond for S} for Theorem 1 is actually setting $\delta$ to be 1 in \eqref{equivalent condition to max type}. Comparing conditions in Theorems \ref{th1}, \ref{th2} and \ref{th3}, it is not hard to see that the union of conditions for $Z_w^\P$ and $Z_\d^\P$ separately is less stringent than those in Theorem \ref{th1}. This is reasonable as Theorem \ref{th1} needs the asymptotic normality of the joint distribution of  $(Z_w^\P,Z_\d^\P)$ while Theorem \ref{th2} and \ref{th3} only needs that for one of the marginal distributions.

For the max-type edge-count test, the limiting distribution of test statistic  still requires the conditions in Theorem \ref{th1}. However, if some techniques are used to conservatively estimate the $p$-value for the max-type statistic, such as the Bonferroni correction, then the union of the conditions in Theorems \ref{th2} and \ref{th3} would be enough.

\begin{corollary} \label{corollary 2}
For graphs with $|G| =O(KN)$ and $1\precsim K \prec N$, if $\var(Q_G)\prec K^{1.5}\sqrt{N}$ and $\ep(Q_G^3) = o(K^2N)$, then, in the usual limiting regime, 
$ Z_w^\P\xrightarrow{\mathcal{D}} N(0,1)$ under the permutation null distribution.
\end{corollary}

\begin{remark}
\textnormal{For the $K$-MST {and the $K$-NNG} on the multivariate data with a fixed dimension, the maximum degree has the order of $K$. Then, $\var(Q_G)$ has the order at most $K^2/4$, and $\ep(Q_G^3)$ has the order at most $K^3$. Thus, the  conditions in Corollary \ref{corollary 2} hold for all $K\prec N$.}

\end{remark}

\begin{corollary} \label{corollary 3}
For graphs with $|G| =O(KN)$ and $1\precsim K \precsim N$, if $1\precsim \var(Q_G)$ and the concentration inequality
\begin{align*}
    P(|Q_G-\ep(Q_G)|\geq t)\leq 2\mathrm{exp}\left(-\frac{ct^2}{{V_G}^a}\right), \quad t>0,
\end{align*}
holds for all large $N$ with some constants $c>0$ and $0<a<1$,
then,
in the usual limit regime, $  Z_\d^\P\xrightarrow{\mathcal{D}} N(0,1) $.
\end{corollary}

\begin{remark}
\textnormal{For the $K$-MST {and the $K$-NNG} constructed on multivariate data with a fixed dimension, conditions in Corollary \ref{corollary 3} hold if $K = O(N^\beta)$, $\beta < 0.5$, because $V_G$ has the order at least $N$ when $1\precsim \var(Q_G)$.}
\end{remark}

\subsection{Original edge-count test}
Consider the usual limit region where $\lim_{N \rightarrow \infty} m/N = r \in (0,1)$. When $r = 0.5$, the original edge-count test is equivalent to the weighted edge-count test asymptotically. Thus, we here study the asymptotic distribution of the original edge-count test statistic when $r \neq 0.5$.

\begin{theorem}\label{th4}In the usual limit region and $\lim_{N\rightarrow \infty} \frac{m}{N} = r$ with $r$ a constant and $r\neq \frac{1}{2}$, under {Condition \ref{cond for Zo}}, $Z_o^\P\xrightarrow{\mathcal{D}} N(0,1)$ under the permutation null distribution.
\end{theorem}

The condition \ref{cond for Zo} in Theorem \ref{th4} are weaker than the conditions required in Theorem \ref{th1} as $T$ has the order of $\max\{|G|,V_G\}$.
\begin{corollary} \label{corollary 4}
For graphs with $|G| =O(KN)$ and $1\precsim K \prec N$, in the usual limit regime and $\lim_{N\rightarrow \infty} \frac{m}{N} = r$ with $r$ a constant and $r\neq \frac{1}{2}$, if either of the following conditions
\begin{itemize}
    \item $|G|$ and $V_G$ do not have the same order, and 
    \begin{align*}
    P(|Q_G-\ep(Q_G)|\geq t)\leq 2\mathrm{exp}\left(-\frac{ct^2}{{T}^a}\right), \quad t>0,
    \end{align*}
    \item $|G| \asymp V_G$, and 
    \begin{align*}
    P(|Q_G-\ep(Q_G)|\geq t)\leq 2\mathrm{exp}\left(-\frac{ct^2}{{N}^a}\right), \quad t>0,
    \end{align*}

\end{itemize}
holds for all large $N$ with some constants $c>0$ and $0<a<1$,
then $  Z_o^\P\xrightarrow{\mathcal{D}} N(0,1) $.
\end{corollary}

\begin{remark}
\textnormal{For the $K$-MST {and the $K$-NNG} constructed on multivariate data with a fixed dimension, conditions in Corollary \ref{corollary 4} hold if $K= O(N^\beta)$ with $\beta < 0.5$ when $|G| \asymp V_G$. When $|G|$ and $V_G$ do not have the same order, these conditions hold as long as $K\prec N$.}
\end{remark}
\subsection{Some brief comments on the conditions}
{
The sufficient conditions in Theorem \ref{th1} are derived  using  Stein's method with the first neighbor dependency.  One key step is to have  the upper bound in the Stein's inequality
\begin{align}\label{eq:bound}
     \sqrt{\frac{2}{\pi}}\ep_\B\bigg|\sum_{i\in \mathcal{N}}\big\{\xi_i\eta_i-\ep_\B(\xi_i\eta_i)\big\} +\sum_{e\in G}\big\{\xi_e\eta_e-\ep_\B(\xi_e\eta_e)\big\} \bigg| + \sum_{i\in \mathcal{N}}\ep_\B|\xi_i\eta_i^2| + \sum_{e\in G}\ep_\B|\xi_e\eta_e^2|
\end{align}
to go to zero, where $\ep_\B$, $\xi_i$, $\eta_i$, $\xi_e$ and $\eta_e$ are defined in Section \ref{proof of th1}. Condition $\sum_{i=1}^N \left|\Tilde{d_i}\right|^3 = o(V_G^{1.5})$ ensures that the quantity $\sum_{i\in \mathcal{N}}\ep_\B|\xi_i\eta_i^2|$ goes to zero. Conditions $\sum_{i=1}^N \Tilde{d_i}^3 = o(V_G\sqrt{|G|})$ and $\sum_{i=1}^N|G_i|^2 = o(|G|^{1.5})$ lead to the zero limit of the quantity $\sum_{e\in G}\ep_\B|\xi_e\eta_e^2|$. To ensure the first quantity in \eqref{eq:bound} to go to zero, under the three previous mentioned conditions ($\sum_{i=1}^N \left|\Tilde{d_i}\right|^3 = o(V_G^{1.5})$, $\sum_{i=1}^N \Tilde{d_i}^3 = o(V_G\sqrt{|G|})$ and $\sum_{i=1}^N|G_i|^2 = o(|G|^{1.5})$), we need two additional conditions,  $\sum_{i=1}^N\sum_{j,k\in node_{G_i}}^{j\neq k}\Tilde{d_j}\Tilde{d_k} = o(|G|V_G)$ and $N_{sq} = o(|G|^2)$.  Conditions in Theorems \ref{th2} and \ref{th4} are derived in a similar way.
}
%{
%The sufficient conditions in Theorem \ref{th2} and \ref{th4} are derived in a similar framework. For the conditions in Theorem \ref{th2}, the quantity $\sum_{i\in \mathcal{N}}\ep_\B|\xi_i\eta_i^2|$ would go to zero unconditionally. The condition $\sum_{i=1}^N|G_i|^2 = o(|G|^{1.5})$ ensures that the quantity $\sum_{e\in G}\ep_\B|\xi_e\eta_e^2|$ would tend to zero, and conditions $\sum_{i=1}^N|G_i|^2 = o(|G|^{1.5})$ and $N_{sq} = o(|G|^2)$ make the first quantity go to zero. For the conditions in Theorem \ref{th4}, conditions $\sum_{i=1}^N|G_i|^2 = o(T^{1.5})$ and $\sum_{i=1}^N\left|\Tilde{d_i}\right|^3 = o(T^{1.5})$ guarantee the quantites $\sum_{i\in \mathcal{N}}\ep_\B|\xi_i\eta_i^2|$ and $\sum_{e\in G}\ep_\B|\xi_e\eta_e^2|$ to go to zero. Then, we deal with the first quantity under these two conditions and derive the condition $\sum_{i=1}^N\sum_{j,k\in node_{G_i}}^{j\neq k}\Tilde{d_j}\Tilde{d_k} = o(T^2)$.}

{
Theorem \ref{th3} is derived in a different way as the test statistic in this case can be expressed into a weighted sum of independent random variables. Then the Lyapunov CLT can be applied.}

The locSCB approach used in proving Theorems \ref{th1}, \ref{th2}, \ref{th4} will be detailed in Section \ref{proof of th1}.  For the moment-based method, it was first proposed in \cite{daniels1944relation}.  Later, \cite{friedman1983graph} claimed that Daniels' conditions can be weakened and provided their new conditions. However, they did not give an explicit proof. \cite{pham1989asymptotic} found that conditions in \cite{friedman1983graph} are not sufficient, so they fixed this problem and proposed a new set of weaker conditions. By using the conditions in \cite{pham1989asymptotic}, the conditions to ensure the asymptotic distributions of  the graph-based tests can be weakened.  However, we found that our locSCB approach  could result in even weaker conditions. A discussion of the conditions from the moment-based method is detailed in Section S2 of the Supplementary Material \citep{zhu2023}.

% \cite{friedman1979} used the moment-based method proposed in \cite{daniels1944relation}. Since then, the moment-based method has been extended. \cite{friedman1983graph} claimed that Daniels' conditions for asymptotic normality can be weakened and gave their relaxed conditions. However, they didn't give an explicit proof and merely stated that they could be obtained by a similar argument as in \cite{daniels1944relation}. \cite{pham1989asymptotic} found that conditions in \cite{friedman1983graph} are not sufficient. They fixed this problem and proposed a new set of relaxed conditions. By using these conditions, the results for the graph-based tests can be relaxed. However, we found that, our `locSCB' approach  could result in much more relaxed conditions for $S$. In the main context in this section, we focus on the `locSCB' approach to prove the Theorem \ref{th1}. The proof of Theorems \ref{th2}, \ref{th3} and \ref{th4} are similar to the proof of Theorem \ref{th1}, and are deferred to Supplement S1. Some comparisons to the moment-based method are provided in the Supplement S2.

% For $Z_o^\P$, we also state the sufficient conditions obtained from the moment-based method in \cite{pham1989asymptotic}.
% the conditions from the `locSCB' approach do not cover these from the moment-based method, so we state another set of conditions from the moment-based method in Theorem \ref{th5}.

\section{Discussions on the new conditions}\label{some discussion on theorems}

\subsection{A comparison to the best existing conditions} \label{comparion of condition}
For the asymptotic distribution of the generalized edge-count test statistic, \cite{chu2019asymptotic} had the best result that required the following conditions:
\begin{align*}
    &|G| = O(N^\alpha)\, 1\leq\alpha<1.25, \quad \sum_{e\in G}|A_e||B_e| = o(|G|^{1.5}), \\ & \sum_{e\in G}|A_e|^2 = o(|G|\sqrt{N}), \quad V_G = O(\sum_{i=1}^N|G_i|^2).
\end{align*}
We here compare our conditions in Theorem \ref{th1} with those in \cite{chu2019asymptotic}. {We first state some propositions with proofs deferred to Section S4 of the Supplementary Material \citep{zhu2023}.}

\begin{enumerate}[label=\textbf{P.\arabic*}]
\item $\sum_{e\in G} |A_e|^2 \asymp \sum_{i=1}^N|G_i|^3$; \label{Prop1}
\item $\sum_{e\in G}|A_e||B_e| \asymp \sum_{i=1}^N|G_i|^2|G_{i,2}|+\sum_{i=1}^N|G_{i,2}|^2$; \label{fact1}
\item $\sum_{i=1}^N\sum_{j,k\in node_{G_i}}^{j\neq k}\Tilde{d_j}\Tilde{d_k} \precsim \sum_{i=1}^N |G_{i,2}|^2+\frac{|G|^2}{N^2}\sum_{i=1}^N|G_i|^2$; \label{fact_djdk}
\item \label{Nsq} $N_{sq}\precsim \sum_{i=1}^N|G_i||G_{i,2}|$.
\end{enumerate}

The condition $\sum_{i=1}^N|G_i|^2 = o(|G|^{1.5})$ in Theorem \ref{th1} can be easily obtained from {Proposition \ref{Prop1}} and under the condition $\sum_{e\in G}|A_e||B_e| = o(|G|^{1.5})$ as $$\sum_{i=1}^N|G_i|^2 \leq \sum_{i=1}^N|G_i|^3 \asymp\sum_{e\in G}|A_e|^2 \precsim \sum_{e\in G}|A_e||B_e|.$$ 
For  conditions $\sum_{i=1}^N|\Tilde{d_i}|^3 = o(V_G^{1.5})$ and $\sum_{i=1}^N \Tilde{d_i}^3 = o(V_G\sqrt{|G|})$ in Theorem \ref{th1}, we have, 
\begin{align} \label{di triple}
    \sum_{i=1}^N \Tilde{d_i}^3\precsim\sum_{i=1}^N |\Tilde{d_i}|^3 \precsim \sum_{i=1}^N|G_i|^3+\frac{|G|}{N}\sum_{i=1}^N|G_i|^2 +\frac{|G|^3}{N^2}.
\end{align}
Under conditions in \cite{chu2019asymptotic}, the right-hand side in (\ref{di triple}) is dominated by $V_G\sqrt{|G|}$, so it is also dominated by $V_G^{1.5}$.
For the condition on $\sum_{i=1}^N\sum_{j,k\in node_{G_i}}^{j\neq k}\Tilde{d_j}\Tilde{d_k}$, from {Proposition \ref{fact_djdk}}, we have 
\begin{align}
    \sum_{i=1}^N\sum_{j,k\in node_{G_i}}^{j\neq k}\Tilde{d_j}\Tilde{d_k}\precsim \sum_{i=1}^N |G_{i,2}|^2+\frac{|G|^2}{N^2}\sum_{i=1}^N|G_i|^2 \precsim \sum_{i=1}^N |G_{i,2}|^2+ V_G|G|\frac{|G|}{N^2}.
\end{align}
Then, from {Proposition \ref{fact1}}, we have $\sum_{i=1}^N |G_{i,2}|^2 \precsim \sum_{e\in G}|A_e||B_e|$, which is dominated by $|G|^{1.5}$ under \cite{chu2019asymptotic}'s conditions.
With the fact that $|G|^{1.5} \precsim |G|V_G$ and $|G| \precsim N^2$, the conditions in \cite{chu2019asymptotic} imply the condition $\sum_{i=1}^N\sum_{j,k\in node_{G_i}}^{j\neq k}\Tilde{d_j}\Tilde{d_k} = o(|G|V_G)$.

For the condition on $N_{sq}$ in Theorem \ref{th1}, from {Propositions \ref{fact1} and \ref{Nsq}}, we have that
\begin{align*}
    N_{sq} \precsim \sum_{i=1}^N|G_i||G_{i,2}| \precsim \sum_{e\in G} |A_e||B_e|.
\end{align*}
Hence, the condition $N_{sq} = o(|G|^{2})$ is weaker than the condition $\sum_{e\in G}|A_e||B_e| = o(|G|^{1.5})$ in \cite{chu2019asymptotic}.

In the above inequalities, many are significantly loosened that the right-hand side could be much larger than the left-hand side.  For a graph that satisfies the conditions in \cite{chu2019asymptotic}, its maximum size needs to be smaller than the order of $N^{1.25}$. But for our new conditions, the size of graph can be much larger.

For the weighted edge-count test, \cite{chen2018weighted} had the best result so far, and required that $\sum_{e\in G}|A_e||B_e| = o(|G|^{1.5})$, $\sum_{e\in G} |A_e|^2 = o(|G|\sqrt{N})$ and $|G| = o(N^{1.5})$. Our conditions in Theorem \ref{th2}  only require $\sum_{i=1}^N|G_i|^2 = o(|G|^{1.5})$ and $N_{sq}\precsim |G|^2$, and they are much weaker as $\sum_{i=1}^N|G_i|^2 \precsim \sum_{e\in G}|A_e||B_e|$ and $N_{sq} \precsim \sum_{e\in G}|A_e||B_e|$.

Existing works have not studied the conditions for $Z_\d^\P$ directly.

For the original edge-count test, \cite{chen2015graph} had the best result so far. They require $\sum_{e\in G}|A_e||B_e| = o(\min\{N^{1.5},|G|^{1.5}\})$ and $|G| = O(N^\alpha), \alpha<1.125$, which are more stringent than conditions in Theorem \ref{th4}. For the condition $\sum_{i=1}^N|\Tilde{d_i}|^3=o(T^{1.5})$ in Theorem \ref{th4}, we have 
$$\sum_{i=1}^N|\Tilde{d_i}|^3\precsim \sum_{i=1}^N|G_i|^3+\frac{|G|}{N}\sum_{i=1}^N|G_i|^2+\frac{|G|^3}{N^2},$$ where the right-hand side is dominated by $|G|^{1.5}$ from {Proposition \ref{fact1}} and under \cite{chen2015graph}'s conditions.
For the condition $\sum_{i=1}^N |G_i|^2 = o(T^{1.5})$ in Theorem \ref{th4}, we have $\sum_{i=1}^N |G_i|^2 \precsim \sum_{e\in G}|A_e||B_e| \prec |G|^{1.5}$. For the condition on $\sum_{i=1}^N\sum_{j,k\in node_{G_i}}^{j\neq k}\Tilde{d_i}\Tilde{d_j}$, from Propositions \ref{fact1} and \ref{fact_djdk} and under \cite{chen2015graph}'s condition on $\sum_{e\in G}|A_e||B_e|$, we have $$\sum_{i=1}^N\sum_{j,k\in node_{G_i}}^{j\neq k}\Tilde{d_i}\Tilde{d_j}\precsim \sum_{i=1}^{N}\left|G_{i, 2}\right|^{2}+\frac{|G|^{2}}{N^{2}} \sum_{i=1}^{N}\left|G_{i}\right|^{2}\prec |G|^2 \precsim T^2.$$ 
For condition $N_{sq} = o(T_1^2)$, we have $N_{sq}\precsim \sum_{e\in G}|A_e||B_e|$ from {Proposition \ref{fact1} and \ref{Nsq}}. Under \cite{chen2015graph}'s condition, $\sum_{e\in G}|A_e||B_e|$ is dominated by $|G|^{1.5}$, so it is also dominated by $T^2$. 

% In the proof of Theorem \ref{th1}-\ref{th3} and Theorem \ref{th4} Set 1, we utilize Stein's method. Hence, 

\subsection{How far are the new conditions from being necessary?}
\begin{figure}[h]
    \centering
    \includegraphics[width = 0.33\textwidth]{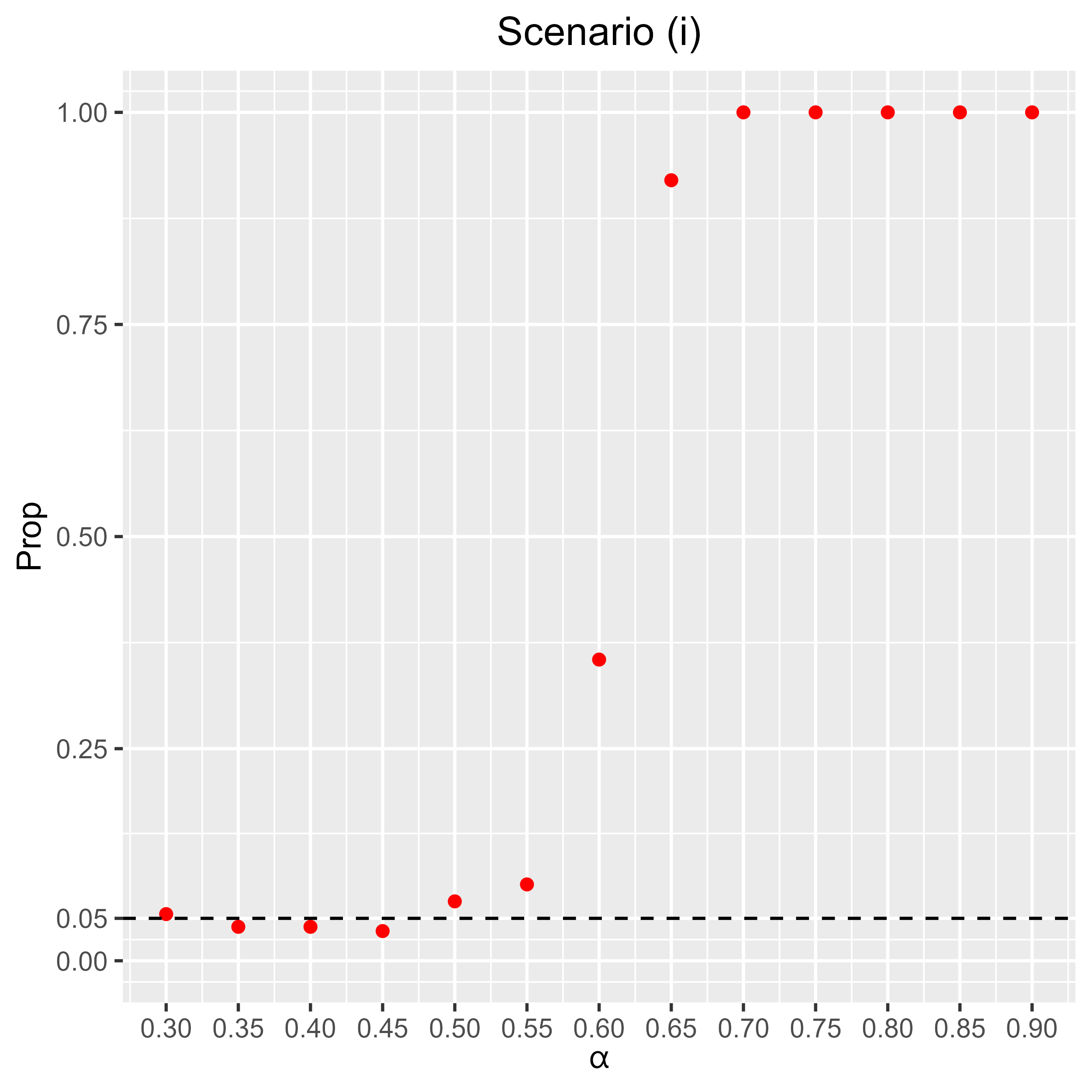}
    \includegraphics[width = 0.33\textwidth]{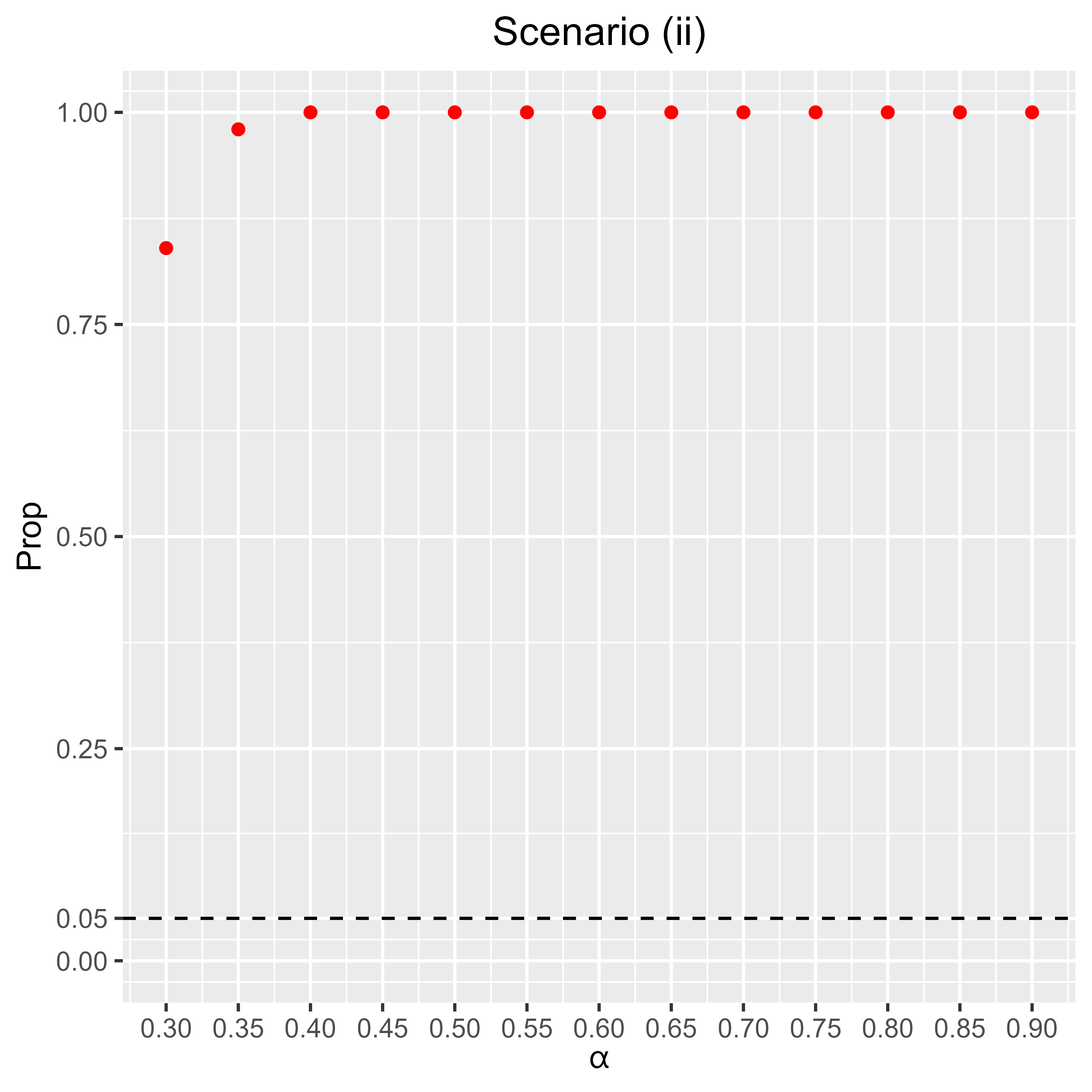}
    \includegraphics[width = 0.33\textwidth]{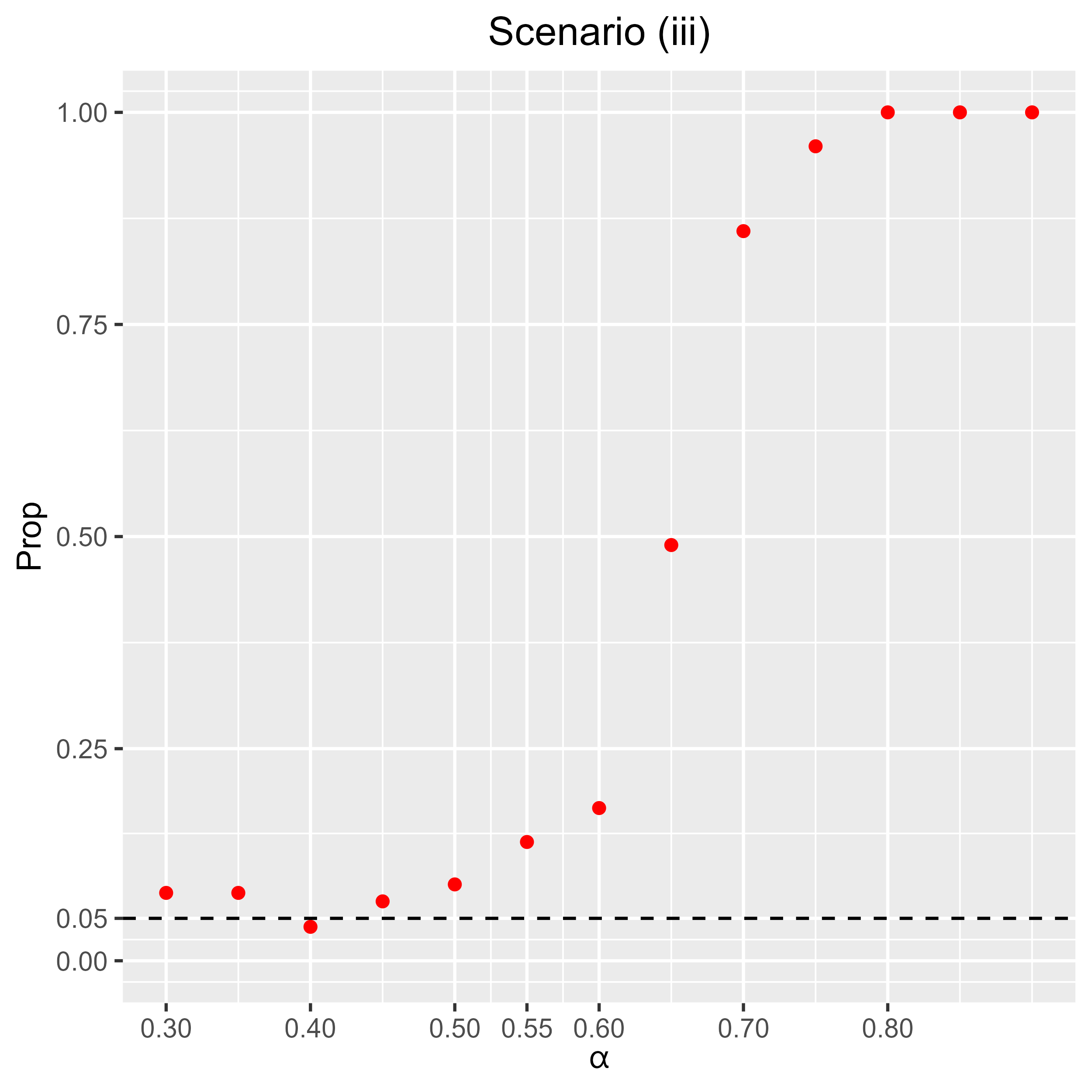}
    \includegraphics[width = 0.33\textwidth]{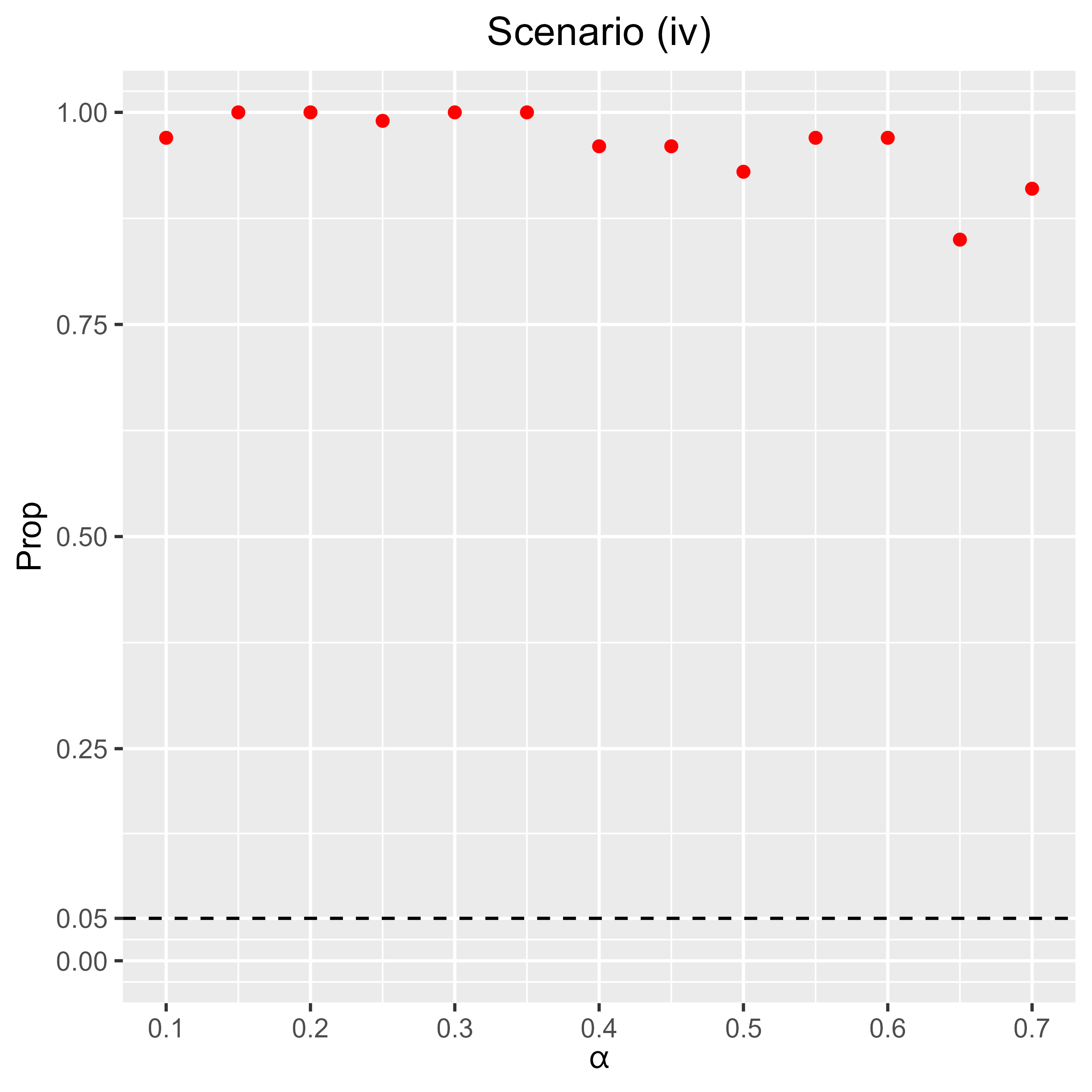}
    \includegraphics[width = 0.33\textwidth]{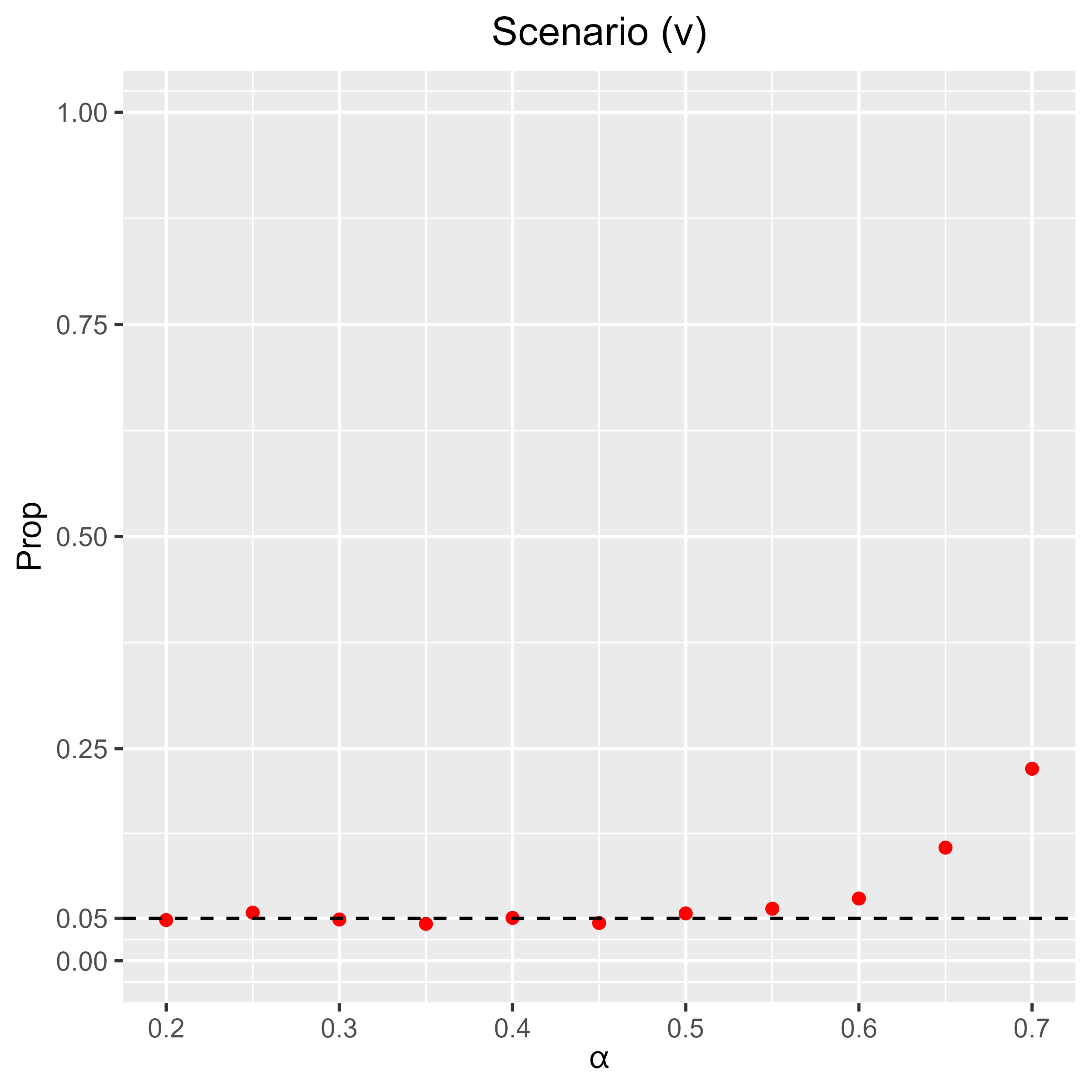}
    \includegraphics[width = 0.33\textwidth]{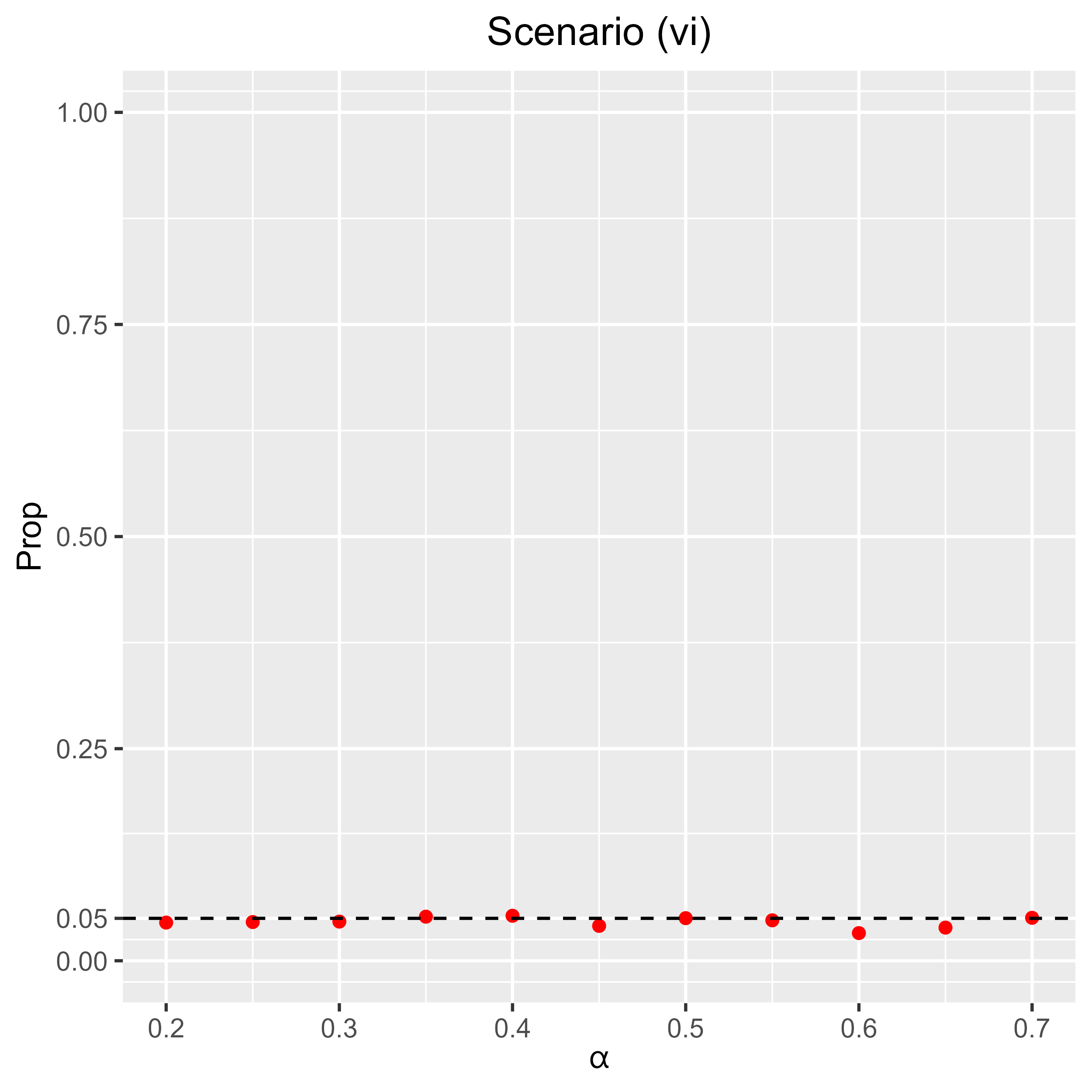}
    \caption{Proportion of rejections from the $\chi^2_2$ distribution under different graph generating rules.}
    \label{check_chi2}
\end{figure}

The conditions in Theorems \ref{th1}, \ref{th2}, \ref{th3} and \ref{th4} are sufficient conditions. One question is how far the new conditions are from being necessary. Here, we focus on the generalized edge-count test and check the validity of the asymptotic $\chi_2^2$ under some synthetic graphs. We repeatedly generate graphs from particular generating rules, obtain the approximate distribution of the test statistic through random permutations, and compare the approximate distribution with the $\chi_2^2$  distribution via Kolmogorov–Smirnov (KS) test.  For each simulation setting, we repeat 100 times, the proportion of rejection of the KS test is plotted in Figure \ref{check_chi2}.
We try to construct graphs that violate the conditions in Theorem \ref{th1}.  The following {six} graph generating rules are considered:
\begin{enumerate}
    \item[(i)] Fix $N$ $(N = 2000)$ nodes indexing from $1, \cdots\, N$. Connect the first node to $N^\alpha$ nodes that are randomly selected from nodes $2, \cdots\, N$. Next, randomly select $N$ edges from all pairwise edges of nodes $2, \cdots\, N$.
    \item[(ii)] Fix $N$ $(N = 2000)$ nodes indexing from $1, \cdots\, N$. Connect the first node to $N^\alpha$ nodes that are randomly selected from nodes $2, \cdots\, N$. Next, connect node $i$ to node $i+1$ for $i\in \{2, \cdots\, N-1\}$ and finally connect nodes $N$ and $2$.
    \item[(iii)] Fix $N$ $(N = 2000)$ nodes indexing from $1, \cdots\, N$. Build the complete graph over the first $M = \lceil N^\alpha \rceil$ nodes. Next, randomly select $2(N-M)$ edges from all pairwise edges of nodes $M+1, \cdots\, N$. Connect the above two subgraphs by adding one edge.
    \item[(iv)] Fix $N$ $(N = 1000)$ nodes indexing from $1,\cdots\, N$, and arrange them to be a circle in the sequence of increasing number. Connect the node $i$ to the next $\lceil N^\alpha \rceil$ nodes. Then connect the node $1$ to node $2+\lceil N^\alpha \rceil$.
    \item[(v)] {Randomly sample $N$ ($N = 2500$) observations from the 2-dimensional standard Gaussian distribution. Build the $K$-MST on the observations with $K = \lceil N^\alpha \rceil$ and $\alpha$ ranging from 0.2 to 0.7.}
    \item[(vi)] {Randomly sample $N$ ($N = 2500$) observations from the 50-dimensional standard Gaussian distribution. Build the $K$-MST on the observations with $K = \lceil N^\alpha \rceil$ and $\alpha$ ranging from 0.2 to 0.7.}
    
\end{enumerate}

Under the graph-generating rule (i), the condition $\sum_{i=1}^N|G_i|^2 = o(|G|^{1.5})$ holds if $0<\alpha<0.75$, the conditions $\sum_{i=1}^N |\tilde{d_i}|^3 = o(V_G^{1.5})$ and $\sum_{i=1}^N \Tilde{d_i}^3 = o(V_G\sqrt{|G|})$ hold if $\alpha<0.5$. The top-left panel in Figure \ref{check_chi2} shows that the $\chi_2^2$ distribution approximation starts to be violated at $\alpha = 0.5$, which is consistent with the analytical result. Under the graph-generating rule (ii), the second sufficient condition $\sum_{i=1}^N |\tilde{d_i}|^3 = o(V_G^{1.5})$ does not hold for any $0<\alpha<1$, which is consistent with the simulation results. Under the graph-generating rule (iii), the condition $\sum_{i=1}^N|G_i|^2 = o(|G|^{1.5})$ holds if $\alpha<0.5$, the condition $\sum_{i=1}^N |\tilde{d_i}|^3 = o(V_G^{1.5})$ holds with $0<\alpha<1$, the condition $\sum_{i=1}^N \Tilde{d_i}^3 = o(V_G\sqrt{|G|})$ holds if $0<\alpha<0.5$ and the condition $N_{sq} = o(|G|^2)$ holds if $\alpha<0.5$, which is consistent with the simulation results. Under the graph-generating rule (iv), the condition $\sum_{i=1}^N |\tilde{d_i}|^3 = o(V_G^{1.5})$ does not hold for any $0<\alpha<1$, which is also consistent with the simulation results in bottom-left panel.

\begin{figure}
    \centering
    \includegraphics[width = 0.75\textwidth]{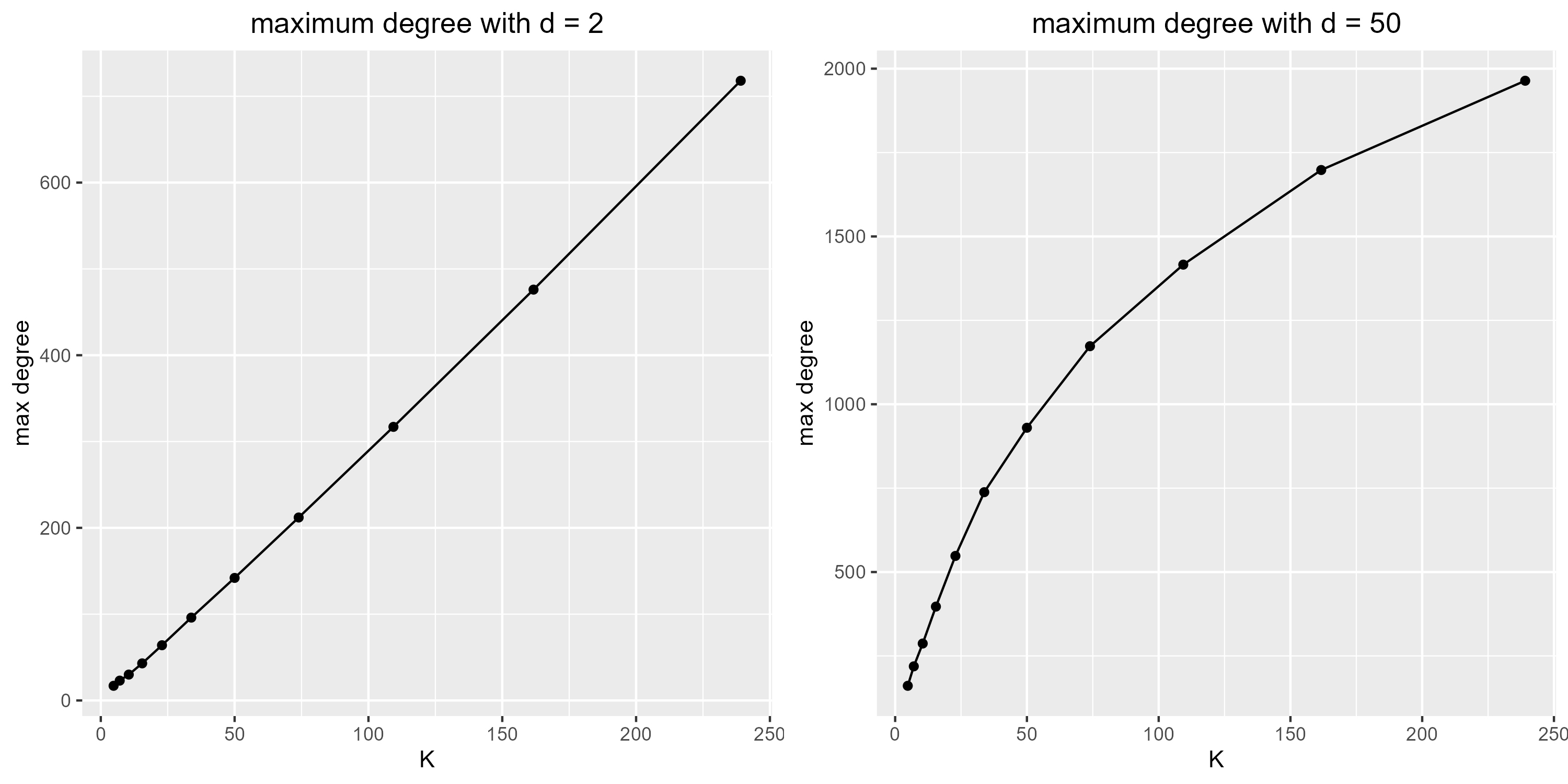}
    \caption{Relationship between the maximum degree of the $K$-MST and $K$ under both dimensions, $d = 2$ (Scenario (v)) and $d=50$ (Scenario (vi)).}
    \label{max_degree_K}
\end{figure}

{Under the graph-generating rules (v) and (vi), the $K$-MST is considered. Remark \ref{remark KMST for S} states that the sufficient conditions in Theorem \ref{th1} hold for the $K$-MST if $K = O(N^\beta)$ with $\beta<0.5$ under fixed dimensions. We can see that the $\chi_2^2$ approximation works well for $\beta < 0.5$. When $d=2$, the approximation starts to deviate when $\alpha$ is bigger than 0.55. Interestingly, when the dimension is larger ($d=50$), the $\chi_2^2$ approximation still works well when $\alpha$ reaches $0.7$. One plausible reason is that Remark \ref{remark KMST for S} considers the asymptotics under fixed dimensions, under which the maximum degree of the $K$-MST is of order $K$.  However, when $d$ is large,  the maximum degree of the $K$-MST is limited to be below the order of $K$ due to the insufficient sample size $N$.  This can be seen from Figure \ref{max_degree_K} where the maximum degree of $K$-MST are plotted over $K$. When $d=50$, if we approximate the maximum degree by $cK^\gamma$, then the estimated value of $\gamma$ is $0.70$. Assume $K = O(N^\beta)$, the conditions in Corollary \ref{corollary 2} hold as long as $2\beta\gamma < 1$ that leads to $\beta < 1/(2\gamma) = 0.713$. This result shows that using the $\chi_2^2$ distribution to approximate the distribution of the generalized edge-count statistic on $K$-MST for data with a relatively large dimension is still a viable option even for $K$'s larger than the upper bound in Remark \ref{remark KMST for S}.}

\begin{figure}
    \centering
    \includegraphics[width = \textwidth]{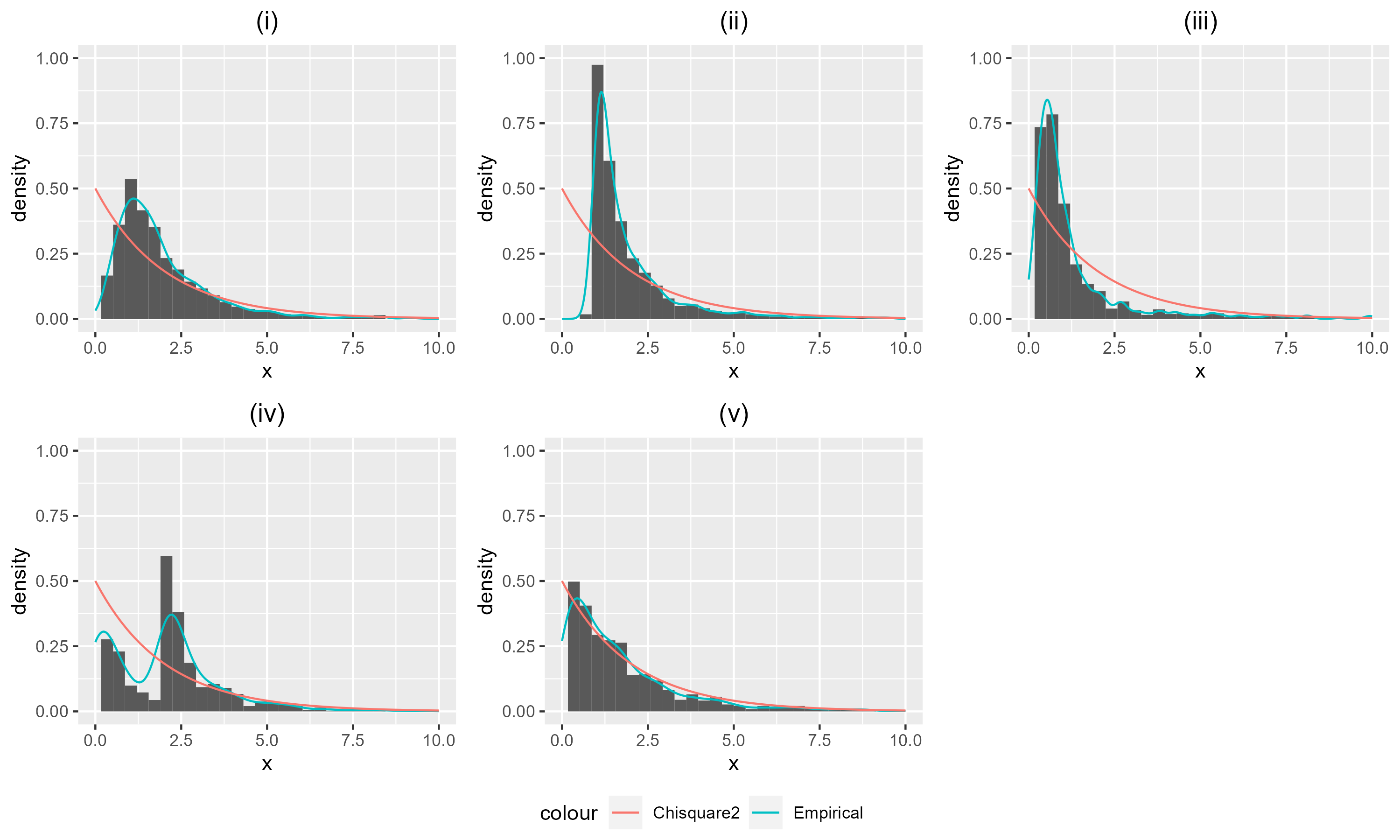}
    \caption{Densities of empirical distributions under scenarios (i) - (v) with specific choices of $\alpha$'s so that the conditions in Theorem \ref{th1} are violated. We set $N = 2000, \alpha=0.8$ for scenarios (i), (ii), (iii),  $N = 1000,  \alpha=0.5$ for scenario (iv), and $N = 2500, \alpha=0.7$ for scenario (v). The red line is the density of the $\chi_2^2$ distribution.}
    \label{density_plot}
\end{figure}

The densities of empirical distributions under these scenarios with specific choices of $\alpha$'s such that the asymptotic $\chi^2_2$ distribution is violated are plotted in Figure \ref{density_plot}. When the asymptotic $\chi^2_2$ distribution is violated, the asymptotic distribution of test statistics under different scenarios can be quite different. However, it seems that the tail probability (on the right) is lighter than the $\chi^2_2$ distribution, which allows us to still use the critical value obtained from the $\chi^2_2$ distribution to control the type I error. Here, the study is through simulations. More systematical investigations will be done in future research.

\subsection{Upper bounds of the difference to the limiting distribution}
\begin{figure}
\includegraphics[width = 0.6\textwidth]{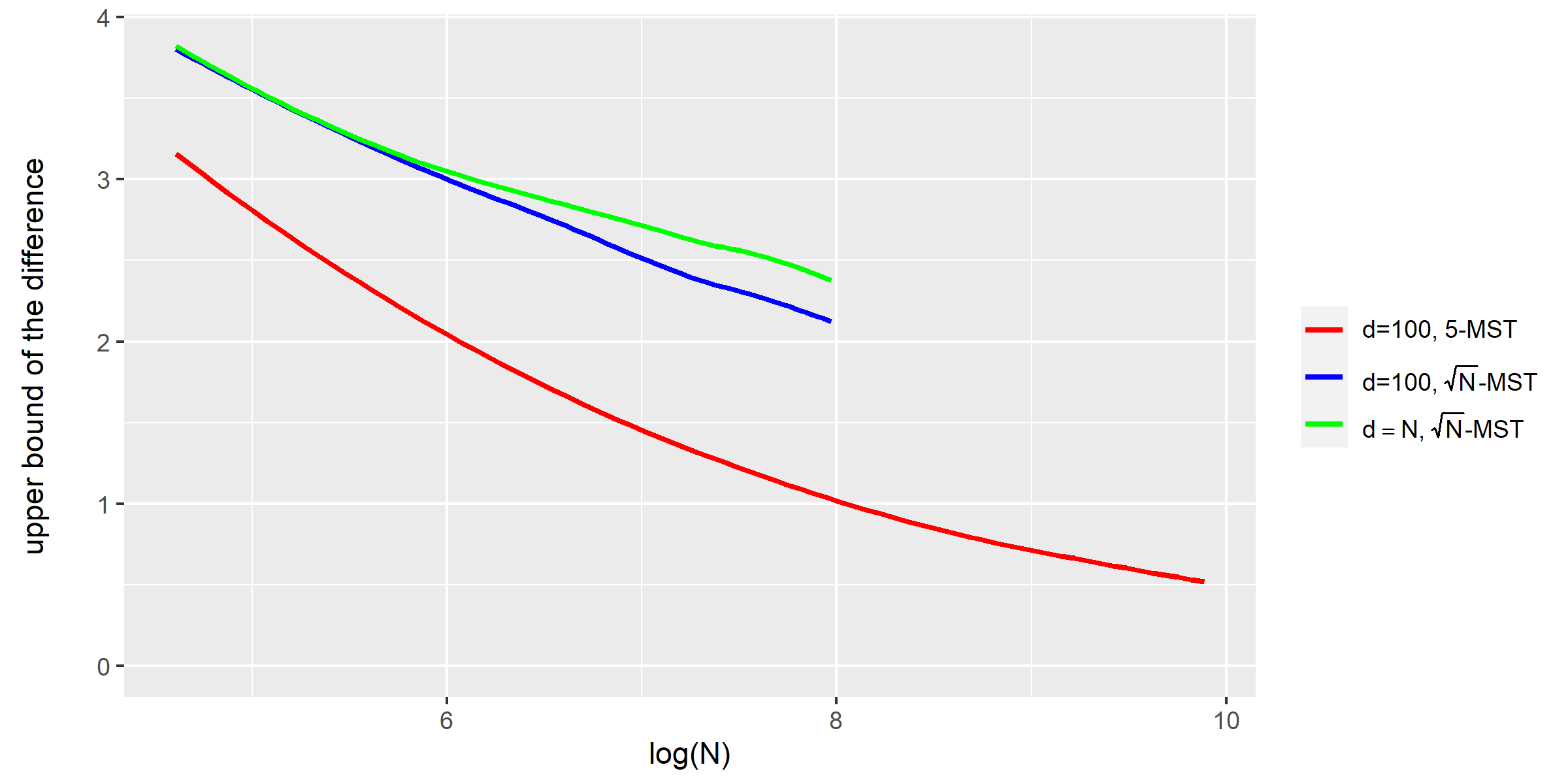}
\caption{\label{upperbound_d100_K5} The upper bound under three settings. Each value in the plot takes the average among 100 simulation runs with the setting. The red line is under the setting (i) with $N$ ranging from 100 to 20000. The blue and green lines are under settings (ii) and (iii), respectively, with $N$ only ranging from 100 to 3000 due to their high computation complexity at large $N$.}
\end{figure}

Since  Stein's method is used, we could compute the upper bound of the difference between the quantity of interest and the standard normal distribution evaluated by Lipschitz-1 functions for finite samples. Here we focus on Theorem \ref{th1}. Figure \ref{upperbound_d100_K5} plots this upper bound (to be more specific, the right-hand side of \eqref{eqn:8} in Section \ref{proof of th1}) for data in different dimensions and graphs in different densities. We consider three settings: (i) $d=100$, $5$-MST; (ii) $d= 100$, $\sqrt{N}$-MST; (iii) $d=N$, $\sqrt{N}$-MST; where  data are all generated from the multivariate Gaussian distribution.  We see that the upper bound decreases as $N$ increases. 

\subsection{{Is a denser graph always more preferable?}}
{Simulation results in Section \ref{merit_denser} show that $\text{GET}_{50}$ has a large power than $\text{GET}_5$ under multiple scenarios. Here, we study the power of the generalized edge-count test on $K$-MST in more detail by increasing the order of $K$ continuously from 0 to 0.85. In particular, we set $m=n=100$ and $K = \lceil N^{\beta} \rceil$ with $\beta$ ranging from 0 to 0.85. We consider 8 different scenarios, among which scenarios (i) - (iv) are the same as those in Section \ref{merit_denser} with a fixed dimension $d = 500$ and scenarios (v) - (ix) are listed below with $\Sigma = (0.5^{|i-j|})_{1\leq i,\space j\leq d}$ and $d = 500$. Scenarios (i) and (ii) compare Gaussian distributions with both mean and variance to be different; scenarios (iii) and (iv) compare non-Gaussian distributions with both mean and variance to be different; scenarios (v)-(vii) further compare Gaussian distributions with only mean difference, only scale difference, and only covariance difference, respectively;  scenario (viii) compares Gaussian and non-Gaussian distributions; and scenario (ix) compares extremely heavy-tailed distributions with both location and scale to be different   .}

\begin{itemize}
    \item[(v) ] $X_1,{\cdots} , X_m \iidsim N(\mathbf{0}_d, \Sigma)$ and $Y_1,{\cdots}, Y_n \iidsim N(0.6/d^{1/3}\mathbf{1}_d, \Sigma)$.
    \item[(vi) ] $X_1,{\cdots} , X_m \iidsim N(\mathbf{0}_d, \Sigma)$ and $Y_1,{\cdots}, Y_n \iidsim N(\mathbf{0}_d, (1+0.17/d^{1/3})^2\times\Sigma)$.
    \item[(vii)] $X_1,{\cdots} , X_m \iidsim N(\mathbf{0}_d, I_d)$ and  $Y_1,{\cdots}, Y_n \iidsim N(\mathbf{0}_d, \Sigma_1)$ where $\Sigma_1 = (0.4^{|i-j|})_{1\leq i,\space j\leq d}$.
    \item[(viii)] $X_1,{\cdots} , X_m \iidsim N(\mathbf{0}_d, I_d)$ and $Y_i = (Y_{i1}^T,Y_{i,2}^T)^T$ with $Y_{i,1} \sim N(\mathbf{0}_{d/2}, I_{d/2})$, $Y_{i,2} \sim t_{30}(\mathbf{0}_{d/2}, I_{d/2})$, and $i=1,{\cdots},n$.
    \item[(ix)] $X_1,{\cdots} , X_m \iidsim t_1(\mathbf{0}_d, I_d)$ and  $Y_1,{\cdots}, Y_n \iidsim t_1(0.6/d^{1/3}\mathbf{1}_d, (1+0.17/d^{1/3})^2\times\Sigma)$.
\end{itemize}

\begin{figure}
    \centering
    \includegraphics[width = 0.6\textwidth]{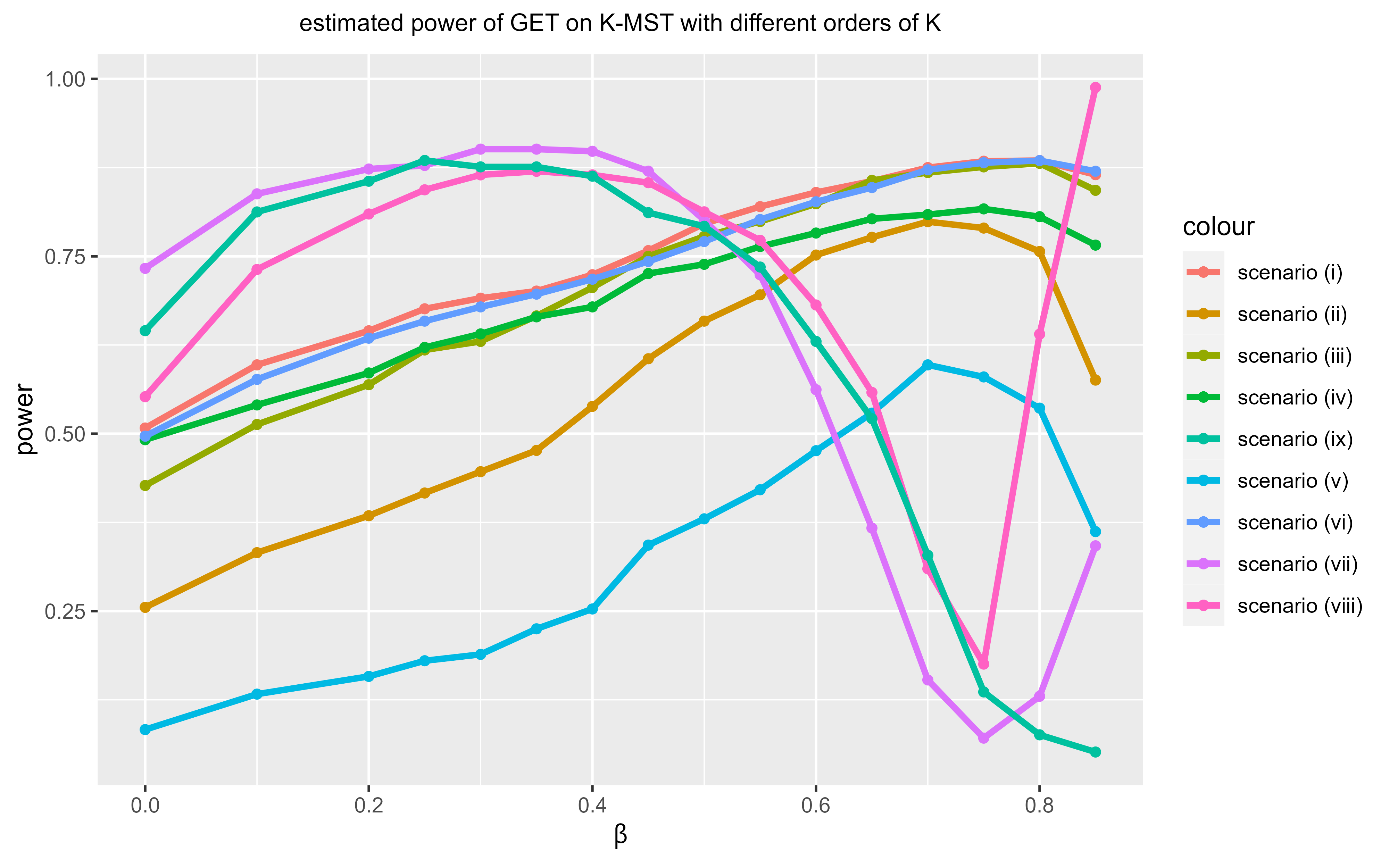}
    \caption{Estimated power at 0.05 significance level of GET on $K$-MST with $K = \lceil N^\beta \rceil$, $\beta\in[0, 0.85]$.}
    \label{power_KMST}
\end{figure}

{For each scenario and each $K$, we run 1,000 trials and the power is estimated as the proportion of trials that reject the null hypothesis at 0.05 significance level.  Figure \ref{power_KMST} plots the estimated power.  We see that, when $\beta$ increases from 0 to 0.25, the power of the test increases for all scenarios.  However, the optimal value of $\beta$ varies across  different scenarios.  For some scenarios, the power increases till $\beta$ reaches 0.8 and then decreases; while for some scenarios, the power begins to decrease at a much smaller $\beta$.  Based on the observation, it is in general safe to consider graphs denser than $O(N)$, while the optimal density of the graph needs further investigation.  One plausible way could be to choose a few representative $K$'s to run the test and then use a multiple testing correction technique, such as the Bonferroni correction, or a $p$-value combining technique, such as the harmonic mean $p$-value, to draw the conclusion. %One plausible and ad-hoc way could be to choose a few representative $K$'s to run the test and draw the conclusion in a holistic way.   % As $\beta$ further increases, the power of the test could decrease for some scenarios Under all scenarios, the estimated power can be improved with properly increasing $K$, but the optimal order of $K$ under different scenarios are not in consistency and it is an open question for future study. It is also worth noticing that the excessive increase in the order of $K$ is detrimental to the power.
}

\subsection{Checking empirical sizes}\label{empirical size section}
Theorems in Section \ref{sec2} provide theoretical guarantees asymptotically. Here, we check empirical sizes for finite samples under a few distributions. Let
\begin{align*}
    &X_i = \Sigma^\frac{1}{2}U_i, i=1, {\cdots}, m, \quad Y_i = \Sigma^\frac{1}{2}V_j,  j =1, {\cdots}, n
\end{align*}
with $\Sigma = (0.5^{|i-j|})_{1\leq i,j\leq d}$. Four distributions for $U_i$'s and $V_j$'s are considered:
\begin{itemize}
    \item[(i)] $U_1, {\cdots}, U_m,V_1,\cdots\,V_n \iidsim N(\mathbf{0}_d,I_d)$,
    \item[(ii)] $U_1, {\cdots}, U_m,V_1,\cdots\,V_n \iidsim t_5(\mathbf{0}_d,I_d)$,
    \item[(iii)] $U_1, {\cdots}, U_m,V_1,\cdots\,V_n \iidsim \text{Uniform} [-1,1]^d$,
    \item[(iv)] $U_1, {\cdots},U_m,V_1,\cdots\,V_n$ are iid with coordinates of $U_i,V_j \iidsim \text{Exp}(1)$.
\end{itemize}
Here, $M =  {N \choose{2}}$ and we use $\lfloor s \rfloor$ to denote the largest integer that is not larger than $s$. 
We consider both the equal sizes $(m=n =100)$ and unequal sizes $(m=50, n=100)$ cases. Table \ref{empirical_GET2} presents the proportion of trials (out of 1,000) that the generalized edge-count test statistic is greater than the $95\%$ quantile of $\chi_2^2$ under the equal sizes case. We see that the empirical size is quite close to the nominal level for all simulation settings. The results for the generalized edge-count test under the unequal sizes case, the original edge-count test, the weighted edge-count test, and the max-type edge-count test are similar, and are presented in Section S4 of the Supplementary Material \citep{zhu2023}.

\begin{table}
\caption{\label{empirical_GET2}Empirical size of the generalized edge-count test with $m = n =100$.}
               \begin{tabular}{|c|c|c|c|c|c|c|c|c|}
        \hline
          distribution&\diagbox{d}{$|G|$}&$\lfloor M^{0.5}\rfloor$&$\lfloor M^{0.6}\rfloor$&$\lfloor M^{0.7}\rfloor$&$\lfloor M^{0.8}\rfloor$&$\lfloor M^{0.9}\rfloor$&$\lfloor M^{0.95}\rfloor$&$\lfloor M^{0.99}\rfloor$ \\
         \hline
        & $0.5N$&0.042& 0.046& 0.055 &0.052 &0.046& 0.054 &0.045\\
         
       (i) normal &$N$&0.049 &0.042& 0.052& 0.056& 0.046& 0.042& 0.059\\
         
         &$1.5N$&0.042 &0.047 &0.037 &0.041 &0.045 &0.056 &0.052\\
         
         &$2N$&0.045& 0.058 &0.049& 0.053& 0.058 &0.051 &0.053\\
         \hline
         & $0.5N$&0.043 &0.042& 0.055 &0.041 &0.052 &0.054& 0.048\\
         
       (ii) $t_5$ &$N$&0.036& 0.04& 0.035 &0.058 &0.055& 0.047 &0.046\\
         
         &$1.5N$&0.034& 0.047 &0.05 &0.06 &0.049 &0.055 &0.047\\
         
         &$2N$&0.045 &0.044 &0.045 &0.053 &0.039& 0.053 &0.066\\
         \hline
        & $0.5N$&0.046& 0.047& 0.044& 0.059& 0.048& 0.046& 0.051\\
         
       (iii) exp(1) &$N$&0.042& 0.047& 0.053 &0.055 &0.05& 0.047& 0.059\\
         
         &$1.5N$&0.047 &0.049 &0.045& 0.045 &0.055 &0.052& 0.052\\
         
         &$2N$&0.05& 0.042 &0.051& 0.052 &0.042 &0.052 &0.051\\
         \hline
         & $0.5N$&0.049 &0.046 &0.048 &0.056 &0.043 &0.048& 0.053\\
         
       (iv) uniform &$N$&0.045& 0.043 &0.042& 0.046& 0.049& 0.047 &0.047\\
         
         &$1.5N$&0.034 &0.056 &0.038& 0.048& 0.048 &0.055 &0.05\\
         
         &$2N$&0.05 &0.047& 0.048 &0.044 &0.055& 0.049 &0.058\\
         \hline
    \end{tabular}
\end{table}

\section{Proof of Theorem \ref{th1}} \label{proof of th1}
%\cite{friedman1979} used the moment-based method proposed in \cite{daniels1944relation}. Since then, the moment-based method has been extended. \cite{friedman1983graph} claimed that Daniels' conditions for asymptotic normality can be weakened and gave their relaxed conditions. However, they didn't give an explicit proof and merely stated that they could be obtained by a similar argument as in \cite{daniels1944relation}. \cite{pham1989asymptotic} found that conditions in \cite{friedman1983graph} are not sufficient. They fixed this problem and proposed a new set of relaxed conditions. By using these conditions, the results for the graph-based tests can be relaxed. However, we found that, our `locSCB' approach  could result in much more relaxed conditions for $S$. In the main context in this section, we focus on the `locSCB' approach to prove the Theorem \ref{th1}. The proof of Theorems \ref{th2}, \ref{th3} and \ref{th4} are similar to the proof of Theorem \ref{th1}, and are deferred to Supplement S1. Some comparisons to the moment-based method are provided in the Supplement S2.

To study the limiting distributions of $Z_w^\P$ and $Z_\d^\P$ jointly, we need to deal with the linear combinations of $\sum_{e\in G}1_{\{J_e=1\}}$ and $\sum_{e\in G}1_{\{J_e=2\}}$. It is clear that the items in these summations are dependent. The dependency comes from two sources. One is due to the permutation null distribution -- given one node from sample X, the probability of another node coming from sample X is no longer $m/N$. The other is due to the nature of the graph-based methods that different edges could share one common node. To conquer these two issues, we work under the bootstrap null distribution to remove the dependency caused by the permutation null distribution first, and then link statistics under the bootstrap null distribution and the permutation null distribution together. For the dependency caused by the nature of the graph-based method, we use the following Stein's method.
\begin{theorem}{(\cite{chen2010normal} Theorem 4.13) }\label{stein}
Let $\{\xi_i,i\in \mathcal{J}\}$ be a random field with mean zero, $W = \sum_{i\in \mathcal{J}}\xi_i$ and $\var(W) = 1$, for each $i\in \mathcal{J}$ there exits $K_i\subset \mathcal{J}$ such that $\xi_i$ and $\xi_{K_i^C}$ are independent, then
$$\sup_{f\in Lip(1)}|\ep f(W)-\ep f(Z)|\leq \sqrt{\frac{2}{\pi}}\ep \bigg|\sum_{i\in \mathcal{J}}\{\xi_i\eta_i-\ep(\xi_i\eta_i)\}\bigg|+\sum_{i\in \mathcal{J}}\ep\bigg|\xi_i\eta_i^2 \bigg|,$$
where $\eta_i = \sum_{j\in K_i}\xi_j$, $Z$ is the standard normal.
\end{theorem}
The Stein's method we rely on is different from that used in \cite{chen2015graph,generalized,chu2019asymptotic}. The main difference is that the theorem used in these earlier works considers a second neighbor of the dependency that, for each $i\in \mathcal{J}$, there exits $K_i\subset L_i\subset \mathcal{J}$ such that $\xi_i$ is independent of $\xi_{K_i^C}$ and $\xi_{K_i}$ is independent of $\xi_{L_i^C}$. Then the upper bound involves $\sum_{j\in L_i}\xi_j$ that could easily expand under the graph structure especially when the graph is dense, causing the conditions to be stringent.  We here turn to the Stein's theorem that only considers the first neighbor of dependency and the resulting quantities can be strategically handled to not expand too much to obtain much weaker conditions.  One difficulty in using this version of Stein's method rather than the second-neighbor version is that the summation is inside the absolute value for the first term on the right-hand side, so we need to manipulate the terms carefully to not relax too much; the detailed manipulations are provided in Section S3 of the Supplementary Material \citep{zhu2023}.

The bootstrap null distribution places probability $1/2^N$ on each of the $2^N$ assignments of $N$ observations to either of the two samples, i.e., each observation is  assigned to sample $X$ with probability $m/N$ and to sample $Y$ with probability $n/N$, independently from any other observations. Let $\ep_\B,\var_\B,\cov_\B$ be expectation, variance and covariance under the bootstrap null distribution. It is  not hard to see that the number of observations assigned to sample $X$ may not be $m$. Let $n_X$ be this number and $Z_X = (n_X-m)/\sigma^\B$ where $\sigma^\B$ is the standard deviation of $n_X$ under the bootstrap null distribution. Notice that the bootstrap null distribution becomes the permutation null distribution conditioning on $n_X=m$. We express $(Z_w^\P,Z_\d^\P)$ in the following way:
\begin{align} \label{eqwd}
    \begin{pmatrix}
    Z_w^\P\\
    Z_\d^\P
    \end{pmatrix}
    =\begin{bmatrix}\sigma_w^\B/\sigma_w^\P &0\\0&\phi_{N}\end{bmatrix}\begin{pmatrix}Z_w^\B\\\sqrt{T_{G}}Z_{\d}^{\B}/\sqrt{V_{G}} \end{pmatrix}+\begin{pmatrix}(\mu_w^\B-\mu_w^\P)/\sigma_w^\P\\ (\mu_\d^\B-\mu_\d^\P)/\sigma_\d^\P\end{pmatrix},
\end{align}
where $Z_w^\B = (R_w-\mu_w^\B)/\sigma_w^\B$, $Z_\d^\B=(R_\d-\mu_\d^\B)/\sigma_\d^\B$, $\phi_{N} = \sqrt{(N-1)/N}$, $T_{G}=\sum_{i=1}^{N}\left|G_{i}\right|^{2}$ and 
\begin{align*}
    &\mu_w^\B = \ep_\B(R_w)=\frac{m n N-m^{2}-n^{2}}{N^{2}(N-2)}|G|,\quad \mu_\d^\B = \ep_\B(R_\d)=\frac{m-n}{N}|G|,\\
    &\sigma_\d^\B = \sqrt{\var_\B(R_\d)} = \sqrt{\frac{mn}{N^2}\sum |G_i|^2}, \\ & \sigma_w^\B = \sqrt{\var_\B(R_w)} = \sqrt{\frac{m^{2} n^{2}}{N^{4}}|G|+\frac{m n}{N^{4}} \frac{(m-n)^{2}}{(N-2)^{2}} \sum_{i=1}^{N}\left|G_{i}\right|^{2}},\\
    &\mu_w^\P= \ep_\P(R_w)=\frac{(n-1)(m-1)}{(N-1)(N-2)}|G|,\quad \mu_\d^\P = \ep_\P(R_\d)=\frac{m-n}{N}|G|,\\
    &\sigma_\d^\P = \sqrt{\frac{m(m-1)n(n-1)}{N(N-1)(N-2)(N-3)}(\frac{m-2}{n-1}+\frac{n-2}{m-1}+2)(\sum |G_i|^2-\frac{4|G|^2}{N})},\\
    &\sigma_w^\P=\sqrt{\frac{m(m-1)n(n-1)}{N(N-1)(N-2)(N-3)}\bigg (|G|-\frac{2}{N(N-1)}|G|^2-\frac{1}{N-2}(\sum |G_i|^2-\frac{4|G|^2}{N})\bigg)}.
\end{align*}
Here we deal with $Z_w^\B$ and $\sqrt{T_G/V_G}Z_\d^\B$ rather than $Z_w^\B$ and $Z_\d^\B$ directly to get rid of the condition $V_G = O(\sum_{i=1}^N|G_i|^2)$ appeared in \cite{generalized} and \cite{chu2019asymptotic}. Since the distribution of $(Z_w^\B,\sqrt{T_{G}/V_{G}} Z_{\d}^{\B})$ under the permutation null distribution is equivalent to the distribution of $(Z_w^\B,\sqrt{T_{G}/V_{G}} Z_{\d}^{\B})|Z_X=0$ under the bootstrap null distribution, we only need to show the following two statements for proving Theorem \ref{th1}:
\begin{enumerate}[label=\textbf{S.\arabic*}]
    \item \label{ result 1} $\Big(Z_{w}^{\B}$, $\sqrt{T_{G}/V_{G}}\left(Z_{\d}^{\B}-\sqrt{1-V_{G}/T_{G}} Z_{X}\right)$, $Z_{X}\Big)$ is asymptotically multivariate Gaussian distributed under the bootstrap null distribution and the covariance matrix of the limiting distribution is of full rank.
    \item \label{ result 2} $\sigma_w^\B/\sigma_w^\P \rightarrow c_w$; $(\mu_w^\B-\mu_w^\P)/\sigma_w^\P \rightarrow 0$ ; $(\mu_\d^\B-\mu_\d^\P)/\sigma_\d^\P\rightarrow 0$ where $c_w$ is a positive constant.
\end{enumerate}
From statement \ref{ result 1}, the asymptotic distribution of $(Z_{w}^{\B}, \sqrt{T_{G}/V_{G}}(Z_{\d}^{\B}-\sqrt{1-V_{G}/T_{G}} Z_{X}))$ conditioning on $Z_X=0$ is a bivariate Gaussian distribution under the bootstrap null distribution, which further implied that the asymptotic distribution of $(Z_w^\B,\sqrt{T_{G}/V_{G}}Z_\d^\B)$ under the permutation null distribution is a bivariate Gaussian distribution. Then, with statement \ref{ result 2}, $\phi_N\rightarrow 1$ and equation (\ref{eqwd}), we have $(Z_w^\P,Z_\d^\P)$ asymptotically bivariate Gaussian distributed under the permutation null distribution. Finally, plus the fact that $\var_\P(Z_w^\P) = \var_\P(Z_\d^\P)=1$ and $\cov_\P(Z_w^\P,Z_\d^\P)=0$, we have that $S\xrightarrow{\mathcal{D}}\chi_2^2$.

The statement \ref{ result 2} is easy to prove. It is clear that $\mu_\d^\B-\mu_\d^\P = 0$. It is also not hard to see that, under condition $\sum_{i=1}^N|G_i|^2 = o(|G|^{1.5})$, in the usual limit regime, $\sigma_w^\B, \sigma_w^\P$ are of the same order of $\sqrt{|G|}$ and $(\mu_w^\B-\mu_w^\P)/\sigma_w^\P$ goes to zero as $|G| = o(N^2)$ and $\sum_{i=1}^N|G_i|^2 = o(|G|N)$.

Next we prove the statement \ref{ result 1}. Let
\begin{align}
    W&=a_{1} Z_{w}^{\B}+a_{2} \sqrt{\frac{T_{G}}{V_{G}}}\left(Z_{\d}^{\B}-\sqrt{1-\frac{V_{G}}{T_{G}}} Z_{X}\right)+a_{3} Z_{X} \nonumber\\
    &=a_{1} Z_{w}^{\B}+a_{2} \sqrt{\frac{T_{G}}{V_{G}}} Z_{\d}^{\B}+\left(a_{3}-a_{2} \sqrt{\frac{T_{G}}{V_{G}}-1}\right) Z_{X}\label{W}.
\end{align}
Firstly we show that, in the usual limit regime,
\begin{eqnarray}
\lim_{N\rightarrow\infty }\var_\B(W)>0 \text{ when at least one of }a_1,a_2,a_3 \text{ is not zero}.
\end{eqnarray}
Since $g_i$'s are independent under the bootstrap null distribution, It is  not hard to derive that
\begin{align*}
    \cov_\B(R_1,n_X) = 2|G|p^2q,\quad \cov_\B(R_2,n_X) = -2|G|pq^2.
\end{align*}
Then we have,
\begin{eqnarray}
&&\cov_\B(Z_w^\B,Z_\d^\B) = pq\frac{(n-m)}{(N-2)}\frac{\sum_{i=1}^N|G_i|^2}{N\sigma_w^\B\sigma_\d^\B}, \nonumber\\
&& \cov_\B(Z_w^\B,Z_X) =  pq\frac{2(n-m)}{(N-2)}\frac{|G|}{N\sigma_w^\B\sigma^\B}, \nonumber\\
&&\cov_\B(Z_\d^\B,Z_X) =  \frac{2pq|G|}{\sigma_\d^\B\sigma^\B}, \quad (\sigma^\B)^2 =Npq,\nonumber
\end{eqnarray}
where $p =m/N, q=n/N$.

Then the variance of $W$ under the bootstrap null distribution can be computed as
\begin{align}
\var_{\B}(W)&=a_{1}^{2}+a_{2}^{2}+a_{3}^{2}+2 a_{1} a_{2} p q \frac{n-m}{N-2} \sqrt{\frac{T_{G}}{V_{G}}} \frac{\sum_{i=1}^{N}\left|G_{i}\right|^{2}}{N \sigma_{w}^{\B} \sigma_{\d}^{\B}} \label{var_W}\\
&\quad +4 p q a_1\frac{n-m}{N-2}\left(a_{3}-a_{2} \sqrt{\frac{T_{G}}{V_{G}}-1}\right) \frac{|G|}{N \sigma_{w}^{\B} \sigma^{\B}}, \nonumber\\
&= a_1^2+a_2^2+a_3^2 + 4a_1a_3\sqrt{pq}\frac{n-m}{N-2}\frac{|G|}{N^{1.5}\sigma_w^\B}+2a_1a_2pq\frac{n-m}{N-2}\frac{\sqrt{V_G}}{\sigma_w^\B N}.\nonumber
\end{align}
Note that $\sigma_w^\B\asymp \sqrt{|G|}$, so $$\frac{|G|}{N^{1.5}\sigma_w^\B} \asymp \frac{\sqrt{|G|}}{N^{1.5}} \rightarrow 0 \text{ and } \frac{\sqrt{V_G}}{\sigma_W^\B N} \precsim \sqrt{\frac{\sum_{i=1}^N |G_i|^2}{|G|N^2}} \rightarrow 0.$$
Thus, we have 
\begin{align*}
    \lim_{N\rightarrow\infty }\var_\B(W) = a_1^2+a_2^2+a_3^2.
\end{align*}
This implies that the covariance matrix of the joint limiting distribution is of full rank. Then by Cram$\acute{\text{e}}$r--Wold device, statement \ref{ result 1} holds if $W$ is asymptotically Gaussian distributed under the bootstrap null distribution for any combinations of $a_1,a_2,a_3$ such that  at least one of them is nonzero.

We reorganized $W$ in the following way
\begin{align*}
    W &= \frac{a_1\Big(\frac{n-1}{N-2}(R_1-|G|p^2)+\frac{m-1}{N-2}(R_2-|G|q^2)\Big )}{\sigma_w^\B}+a_2\sqrt{\frac{T_G}{V_G}}\frac{\Big( R_1-R_2-|G|p^2+|G|q^2\Big)}{\sigma_\d^\B}\\
      &\quad+\Big(a_3-a_2\sqrt{\frac{T_G}{V_G}-1}\Big)\frac{(n_X-m)}{\sigma^\B}\\
      &= \sum_{e\in G}\Big(\frac{a_1}{\sigma_w^\B}\frac{N}{N-2}(I_{\{g_{e^+}=1\}}-p)(I_{\{g_{e^-}=1\}}-p)-\frac{a_1}{\sigma_w^\B}\frac{I_{\{J_e=1\}}+I_{\{J_e=2\}}-p^2-q^2}{N-2}\Big)\\
      &\quad+\sum_{e\in G}a_2\sqrt{\frac{T_G}{V_G}}\frac{1}{\sigma_\d^\B}\Big(I_{\{g_{e^+}=1\}}+I_{\{g_{e^-}=1\}}-2p\Big)+\sum_{i\in \mathcal{N}}\Big(a_3-a_2\sqrt{\frac{T_G}{V_G}-1}\Big)\frac{(I_{\{g_i=1\}}-p)}{\sigma^\B}.
\end{align*}
Define a function $h: \mathcal{N}\rightarrow \mathbb{R}$ and $h(i) = I_{\{g_i=1\}}-p,  i \in \mathcal{N}$. Then,
\begin{eqnarray}
&& (I_{\{g_{e^+}=1\}}-p)(I_{\{g_{e^-}=1\}}-p) = h(e^{+})h(e^-),\nonumber\\
&& I_{\{J_e=1\}}+I_{\{J_e=2\}}-p^2-q^2=2h(e^+)h(e^-)+(p-q)(h(e^+)+h(e^-)),\nonumber\\
&& I_{\{g_{e^+}=1\}}+I_{\{g_{e^-}=1\}}-2p= h(e^+)+h(e^-).\nonumber
\end{eqnarray}
Thus, $W$ can be expressed as
\begin{align*}
    W =& \sum_{e\in G}\frac{a_1}{\sigma_w^\B}h(e^+)h(e^-)+\left(\frac{a_2}{\sigma_\d^\B}\sqrt{\frac{T_G}{V_G}}-\frac{a_1(p-q)}{\sigma_w^\B(N-2)}\right)\sum_{i=1}^N |G_i|h(i)+ \sum_{i\in \mathcal{N}}\Big(a_3-a_2\sqrt{\frac{T_G}{V_G}-1}\Big)\frac{h(i)}{\sigma^\B}\\
=&\sum_{e\in G}\frac{a_1}{\sigma_w^\B}h(e^+)h(e^-)+\sum_{i=1}^N\left(\frac{a_2}{\sqrt{pqV_G}}\left(|G_i|-\frac{2|G|}{N}\right)-\frac{a_1(p-q)|G_i|}{\sigma_w^\B(N-2)}+\frac{a_3}{\sqrt{pqN}}\right)h(i).
\end{align*}

Let $\xi_e = b_0h(e^+)h(e^-)$ and $\xi_i = b_ih(i)$ with
\begin{align*}
    b_0 = \frac{a_1}{\sigma_w^\B}, \quad b_i = \frac{a_2}{\sqrt{pqV_G}}\left(|G_i|-\frac{2|G|}{N}\right)-\frac{a_1(p-q)|G_i|}{\sigma_w^\B(N-2)}+\frac{a_3}{\sqrt{pqN}},
\end{align*}
for $i=1,{\cdots},N$.
Then
\begin{equation}
\label{eqn:W}
    W = \sum_{e\in G}\xi_e+\sum_{i\in \mathcal{N}}\xi_i.
\end{equation}
Plugging in the expressions of $\sigma_w^\B$, It is  not hard to see that
$$|b_0|\precsim \frac{1}{\sqrt{|G|}},\quad |b_i|\precsim \frac{1}{\sqrt{V_G}}\left||G_i|-\frac{2|G|}{N}\right|+\frac{1}{\sqrt{N}}. $$

Next, we apply the Stein's method to $\widetilde{W} = W/\sqrt{\var_\B(W)}$. Let $\mathcal{J} = G\cup \mathcal{N}$, $K_e = A_e\cup \{e^+,e^-\}$ for each edge $e = (e^+,e^-)\in G$ and $K_i =  G_i\cup\{i\}$ for each node $i\in \mathcal{N}$. These $K_e \text{'s}$ and $K_i\text{'s}$ satisfy the assumptions in Theorem \ref{stein} under the bootstrap null distribution. Then, we define $\eta_e\text{'s}$ and $\eta_i\text{'s}$ as follows:
\begin{align*}
    \eta_e = \xi_{e^+}+\xi_{e^-}+\sum_{e\in A_e}\xi_e,\text{ for each edge } e\in G, \text{ and }\eta_i = \xi_i+\sum_{e\in G_i}\xi_e, \text{ for each node } i\in \mathcal{N} .
\end{align*}
By theorem \ref{stein}, we have
\begin{eqnarray}
\label{eqn:8}
\sup_{h\in Lip(1)}|\ep_\B h(\widetilde{W})-\ep_\B h(Z)| &\leq& \sqrt{\frac{2}{\pi}}\frac{1}{\var_\B(W)}\ep_\B\bigg|\sum_{i\in \mathcal{N}}\{\xi_i\eta_i-\ep_\B(\xi_i\eta_i)\} +\sum_{e\in \mathcal{G}}\{\xi_e\eta_e-\ep_\B(\xi_e\eta_e)\} \bigg| \nonumber\\
& &+\frac{1}{\var_\B^{1.5}(W)}\Big(\sum_{i\in \mathcal{N}}\ep_\B|\xi_i\eta_i^2|+\sum_{e\in G}\ep_\B|\xi_e\eta_e^2|\Big).
\end{eqnarray}

Our next goal is to find some conditions under which the RHS\footnote{RHS: right-hand side; LHS: left-hand side.} of inequality (\ref{eqn:8}) can go to zero. Since the limit of $\var_\B(W)$ is bounded above zero when $a_1,a_2,a_3$ are not all zeros, the RHS of inequality (\ref{eqn:8}) goes to zero if the following three terms go to zero:
\begin{enumerate}
    \item[(A1)]  
    $\ep_\B\bigg|\sum_{i\in \mathcal{N}}\big\{\xi_i\eta_i-\ep_\B(\xi_i\eta_i)\big\} +\sum_{e\in G}\big\{\xi_e\eta_e-\ep_\B(\xi_e\eta_e)\big\} \bigg|$,
    \item[(A2)] 
    $\sum_{i\in \mathcal{N}}\ep_\B|\xi_i\eta_i^2| $,
    \item[(A3)]
    $\sum_{e\in G}\ep_\B|\xi_e\eta_e^2|$.
\end{enumerate}

For (A1), we have
\begin{eqnarray}\label{eqn: 9}
&&\ep_\B\bigg|\sum_{i\in \mathcal{N}}\big\{\xi_i\eta_i-\ep_\B(\xi_i\eta_i)\big\} +\sum_{e\in G}\big\{\xi_e\eta_e-\ep_\B(\xi_e\eta_e)\big\} \bigg|\nonumber \\ 
&&\leq \ep_\B\bigg|\sum_{i\in \mathcal{N}}\big\{\xi_i\eta_i-\ep_\B(\xi_i\eta_i)\big\}|\bigg|+\ep_\B\bigg|\sum_{e\in G}\big\{\xi_e\eta_e-\ep_\B(\xi_e\eta_e)\big\} \bigg|\nonumber\\
&&\leq \sqrt{\sum_{i\in \mathcal{N}}\var_\B(\xi_i\eta_i)+\sum_{i,j}^{i\neq j}\cov_\B(\xi_i\eta_i,\xi_j\eta_j)}\nonumber+\sqrt{\sum_{e\in G}\var_\B(\xi_e\eta_e)+\sum_{e,f}^{e\neq f}\cov_\B(\xi_e\eta_e,\xi_f\eta_f)} \nonumber\\
&&= \sqrt{\sum_{i\in \mathcal{N}}\var_\B(\xi_i\eta_i)+\sum_{i=1}^N\sum_{j\in node_{G_{i,2}}}\cov_\B(\xi_i\eta_i,\xi_j\eta_j)}\nonumber\\
&&\hspace{5mm}+\sqrt{\sum_{e\in G}\var_\B(\xi_e\eta_e)+\sum_{e \in G}\sum_{f\in C_e\backslash\{e\}}\cov_\B(\xi_e\eta_e,\xi_f\eta_f)} \nonumber.
\end{eqnarray}

The last equality holds as $\xi_i\eta_i$ and $\{\xi_j\eta_j\}_{j\notin node_{G_{i,2}}}$ are uncorrelated and $\xi_e\eta_e$ and  $\{\xi_f\eta_f\}_{f\notin C_e}$ are uncorrelated under the bootstrap null distribution. The covariance part of edges is a bit complicated to handle directly, so we decompose it into three parts based on the relationship of $e$ and $f$:
\begin{eqnarray*}
\sum_{e \in G}\sum_{f\in C_e\backslash\{e\}}\cov_\B(\xi_e\eta_e,\xi_f\eta_f)&&=\sum_{e \in G}\sum_{f\in A_e\backslash\{e\}}\cov_\B(\xi_e\eta_e,\xi_f\eta_f)
+\sum_{e \in G}\sum_{f\in B_e\backslash A_e}\cov_\B(\xi_e\eta_e,\xi_f\eta_f) \\
&& \quad +\sum_{e \in G}\sum_{f\in C_e\backslash B_e}\cov_\B(\xi_e\eta_e,\xi_f\eta_f). \nonumber
\end{eqnarray*}

By carefully examining these quantities, we can show the following inequalities (\ref{cl2}) - (\ref{cl9}). Here, we only need to consider the worst case, i.e. $b_0$ and $b_i$'s take the largest possible orders, denoted by $c_0$ and $c_i$'s, i.e. $$c_0 = \frac{1}{\sqrt{|G|}},\text{ and }c_i= \frac{1}{\sqrt{V_G}}\left||G_i|-\frac{2|G|}{N}\right|+\frac{1}{\sqrt{N}}.$$ The details of obtaining (\ref{cl2}) - (\ref{cl9}) are provided in Section S3 of the Supplementary Material \citep{zhu2023}.
\begin{eqnarray}
&&\sum_{i=1}^N\varb\left(\xi_i\eta_i\right) \precsim \sum_{i=1}^Nc_i^4+c_0^2\sum_{i=1}^Nc_i^2|G_i|,\label{cl2}\\
&&\sum_{e \in G} \var_{\B}\left(\xi_{e} \eta_{e}\right) \precsim c_0^2\sum_{i=1}^N c_i^2|G_i| + c_0^3\sum_{i=1}^N c_i|G_i| + c_0^4 \sum_{i=1}^N|G_i|^2\label{cl3}\\
&&\sum_{i=1}^N\sum_{j\in node_{G_{i,2}}} \cov_\B(\xi_i\eta_i,\xi_j\eta_j) \precsim \sum_{i=1}^N\sum_{j\in node_{G_i}}\big(c_0c_i^2c_j+c_0^2c_ic_j)\label{cl4}\\
&&\hspace{4.5cm}+\left|c_0^2 \sum_{i=1}^N\sum^{j\in node_{G_{i,2}}}b_ib_jN_{i,j}\right|\nonumber\\
&&\sum_{e\in G}\sum_{f\in A_e\backslash\{e\}}\cov_\B(\xi_e\eta_e,\xi_f\eta_f)\precsim\sum_{i=1}^N\sum_{j,k\in node_{G_i}}^{j\neq k}\left(c_0^3c_j+c_0^3c_k+c_0^3c_i1_{\{(j,k) \text{ exists}\}}\right)\label{cl5}\\
&&\hspace{5cm}+c_0^4N_{sq}+\left|\sum_{i=1}^N\sum_{j,k\in node_{G_i}}^{j\neq k}c_0^2b_jb_k\right|\nonumber\\
&&\sum_{e\in G}\sum_{f\in B_e\backslash A_e}\cov_\B(\xi_e\eta_e,\xi_f\eta_f)\precsim  c_0^4N_{sq}\label{cl6}\\
&&\sum_{e\in G}\sum_{f\in C_e\backslash B_e} \cov_\B(\xi_e\eta_e,\xi_f\eta_f) = 0,\label{cl7}\\
&&\sum_{i=1}^N \ep_\B[|\xi_i\eta_i^2|] \precsim \sum_{i=1}^Nc_i^3+c_0^2\sum_{i=1}^Nc_i|G_i|,\label{cl8}\\
&&\sum_{e\in G}\ep_\B[|\xi_e|\eta_e^2] \precsim  |G|c_0^3+ c_0\sum_{i=1}^N |G_i|c_i^2+c_0^3\sum_{i=1}^N|G_i|^2\label{cl9}.
\end{eqnarray}
Based on facts that $c_i\precsim 1$ for all $i$'s, (A1), (A2) and (A3) going to zero as long as the following conditions hold.

\begin{multicols}{2}
\noindent
    \begin{equation}
        \sum_{i=1}^N c_i^3 \rightarrow 0,\label{si0}
    \end{equation}
    \begin{equation}
        c_0^2\sum_{i=1}^Nc_i|G_i| \rightarrow 0, \label{si1}
    \end{equation}
    \begin{equation}
        |G|c_0^3 \rightarrow 0,\label{si2}
    \end{equation}
    \begin{equation}
        c_0\sum_{i=1}^N|G_i|c_i^2\rightarrow 0,\label{si3}
    \end{equation}
    \begin{equation}
        c_0^3\sum_{i=1}^N|G_i|^2 \rightarrow 0,\label{si4}
    \end{equation}
    \begin{equation}
        \sum_{i=1}^N\sum_{j\in node_{G_i}}\big(c_0c_i^2c_j+c_0^2c_ic_j)\rightarrow 0,\label{si5}
    \end{equation}
    \begin{equation}
        c_0^2 \sum_{i=1}^N\sum_{j\in node_{G_{i,2}}}b_ib_jN_{i,j} \rightarrow 0,\label{si6}
    \end{equation}
    \begin{equation}
        c_0^2\sum_{i=1}^N\sum_{j,k\in node_{G_i}}^{j\neq k}b_ib_j \rightarrow 0, \label{si6-2}
    \end{equation}
        \begin{equation}
        c_0^4(N_{sq}+\sum_{i=1}^N|G_i|^2)\rightarrow 0,\label{si8}
    \end{equation}
\end{multicols}
\begin{equation}
        \sum_{i=1}^N\sum_{j,k\in node_{G_i}}^{j\neq k}\left(c_0^3c_j+c_0^3c_k+c_0^3c_i1_{\{(j,k) \text{ exists}\}}\right) \rightarrow 0.\label{si7}
\end{equation}

%\begin{align}
%    &\sum_{i=1}^N c_i^3 \rightarrow 0\label{si0},\\
%    &c_0^2\sum_{i=1}^Nc_i|G_i| \rightarrow 0 \label{si1},\\
%    &|G|c_0^3 \rightarrow 0\label{si2},\\
%    &c_0\sum_{i=1}^N|G_i|c_i^2\rightarrow 0\label{si3},\\
%    &c_0^3\sum_{i=1}^N|G_i|^2 \rightarrow 0\label{si4},\\
%    &\sum_{i=1}^N\sum_{j\in node_{G_i}}\big(c_0c_i^2c_j+c_0^2c_ic_j)\rightarrow 0\label{si5},\\
%    &c_0^2 \sum_{i=1}^N\sum_{j\in node_{G_{i,2}}}b_ib_jN_{i,j} \rightarrow 0\label{si6},\\
%    &c_0^2\sum_{i=1}^N\sum_{j,k\in node_{G_i}}^{j\neq k}b_ib_j \rightarrow 0 \label{si6-2}\\
%    &\sum_{i=1}^N\sum_{j,k\in node_{G_i}}^{j\neq k}\left(c_0^3c_j+c_0^3c_k+c_0^3c_i1_{\{(j,k) \text{ exists}\}}\right) \rightarrow 0\label{si7},\\
%    &c_0^4(N_{sq}+\sum_{i=1}^N|G_i|^2)\rightarrow 0\label{si8}.
%\end{align}

Next, we show that the conditions in Theorem \ref{th1} can ensure (\ref{si0})-(\ref{si7}).

For  condition (\ref{si0}), we have 
\begin{align*}
    \sum_{i=1}^Nc_i^3 = \sum_{i=1}^N\left(\frac{1}{\sqrt{V_G}}\left||G_i|-\frac{2|G|}{N}\right|+\frac{1}{\sqrt{N}}\right)^3 \precsim \sum_{i=1}^N\frac{||G_i|-\frac{2|G|}{N}|^3}{V_G^{1.5}}+\frac{1}{\sqrt{N}},
\end{align*}
so  condition (\ref{si0}) holds when $\sum_{i=1}^N\left|\Tilde{d}_i\right|^3/V_G^{1.5}\rightarrow 0$.

For  condition (\ref{si1}), we have
\begin{align*}
    c_0^2\sum_{i=1}^Nc_i|G_i|= c_0^2\sum_{i=1}^N|G_i|\frac{1}{\sqrt{V_G}}\left||G_i|-\frac{2|G|}{N}\right|+c_0^2\frac{1}{\sqrt{N}}\sum_{i=1}^N|G_i|.
\end{align*}
It is easy to check that the second term has the order of $1/\sqrt{N}$, so it goes to zero as $N \rightarrow \infty$. The first term is dominated by $\max|\Tilde{d}_i|/\sqrt{V_G}$, which goes to zero under  condition $\sum_{i=1}^N\left|\Tilde{d}_i\right|^3/V_G^{1.5}\rightarrow 0$ from Theorem 1 in \cite{hoeffding1951combinatorial} with $r$ taking $3$. Thus,  condition (\ref{si1}) holds when $\sum_{i=1}^N\left|\Tilde{d}_i\right|^3/V_G^{1.5}\rightarrow 0$.

Condition (\ref{si2}) holds trivially as $|G|c_0^3 = 1/\sqrt{|G|}$.

For  condition (\ref{si3}), we have 
\begin{align*}
 c_0\sum_{i=1}^N|G_i|c_i^2 = c_0\sum_{i=1}^N|G_i|\left(\frac{1}{\sqrt{V_G}}\left||G_i|-\frac{2|G|}{N}\right|+\frac{1}{\sqrt{N}}\right)^2  \precsim c_0\sum_{i=1}^N|G_i|\frac{\left(|G_i|-\frac{2|G|}{N}\right)^2}{V_G}+c_0\frac{1}{N}\sum_{i=1}^N|G_i|.
\end{align*}
The second term has the order of $\sqrt{|G|}/N$, so it goes to zero when $|G| = o(N^2)$. The first term can be written as 
\begin{align*}
    c_0\sum_{i=1}^N|G_i|\frac{\left(|G_i|-\frac{2|G|}{N}\right)^2}{V_G} = c_0\frac{\sum_{i=1}^N\left(|G_i|-\frac{2|G|}{N}\right)^3}{V_G}+c_0\frac{2|G|}{N},
\end{align*}
so it goes to zero under conditions $\sum_{i=1}^N\Tilde{d}_i^3/(V_G\sqrt{|G|})\rightarrow 0 $ and $|G| = o(N^2)$.

For  condition (\ref{si4}), it is directly equivalent to $\sum_{i=1}^N|G_i|^2 = o(|G|^{1.5})$. Note that  condition $|G| = o(N^2)$ also holds under  condition $\sum_{i=1}^N|G_i|^2 = o(|G|^{1.5})$ as $\sum_{i=1}^N|G_i|^2\geq 4|G|^2/N$.

For  condition (\ref{si5}), we have with the fact $c_i\precsim 1$ for all $i$'s 
\begin{align*}
    &\sum_{i=1}^N\sum_{j\in node_{G_i}}c_0c_i^2c_j \precsim \sum_{i=1}^N\sum_{j\in node_{G_i}}c_0c_i^2 = c_0\sum_{i=1}^N|G_i|c_i^2,\\
    &\sum_{i=1}^N\sum_{j\in node_{G_i}}c_0^2c_ic_j \precsim \sum_{i=1}^N\sum_{j\in node_{G_i}}c_0^2c_i = c_0^2\sum_{i=1}^N|G_i|c_i,
\end{align*}
where both the right-hand sides go to zero from (\ref{si1}) and (\ref{si3}).

For  condition (\ref{si6}), we have 
\begin{align*}
    \sum_{i=1}^N\sum_{j\in node_{G_{i,2}}}b_ib_jN_{ij} &= \sum_{i=1}^N\sum_{j\in node_{G_{i,2}}}b_ib_j\sum_{k=1}^N1_{\{k \text{ connects i, j simultaneously }\}}\\
    &= \sum_{k=1}^N\sum_{i\in node_{G_k}}\sum_{j\in node_{G_k}\backslash\{i\}}b_ib_j = \sum_{k=1}^N\sum_{i,j\in node_{G_k}}^{i\neq j}b_ib_j,
\end{align*}
so the left-hand side of  condition (\ref{si6}) is equal to the term in  condition (\ref{si6-2}).

For  condition (\ref{si6-2}), we have 
\begin{align*}
     c_0^2\sum_{i=1}^N\sum_{j,k\in node_{G_i}}^{j\neq k}b_jb_k \precsim& \left|\frac{1}{|G|V_G} \sum_{i=1}^N\sum_{j,k\in node_{G_i}}^{j\neq k}\left(|G_j|-\frac{2|G|}{N}\right)\left(|G_k| - \frac{2|G|}{N}\right)\right|\\
     &+\frac{\sum_{i=1}^N|G_i|^{1.5}}{|G|\sqrt{N}} +\frac{\sum_{i=1}^N|G_i|^2}{|G|N}.
\end{align*}
It is  not hard to see that the second and the third term would go to zero as $N\rightarrow \infty$ under  condition $\sum |G_i|^2/|G|^{1.5} \rightarrow 0$ as
\begin{align*}
    \frac{\sum_{i=1}^N|G_i|^{1.5}}{|G|\sqrt{N}}\precsim \frac{\sqrt{|G|\sum_{i=1}^N|G_i|^2}}{|G|N},\quad \frac{\sum_{i=1}^N|G_i|^2}{|G|N} \precsim \frac{\sum_{i=1}^N|G_i|^2}{|G|^{1.5}}\frac{\sqrt{|G|}}{N}.
\end{align*}
Thus  condition (\ref{si6-2}) holds if  conditions $\sum_{i=1}^N\sum_{j,k\in node_{G_i}}^{j\neq k}\left(|G_j|-\frac{2|G|}{N}\right)\left(|G_k| - \frac{2|G|}{N}\right)$ $= o(|G|V_G)$ and $\sum_{i=1}^N|G_i|^2= o(|G|^{1.5})$ hold.

For  condition (\ref{si8}), it directly requires $N_{sq} = o(|G^2|)$ as $c_0^3\sum_{i=1}^N|G_i|^2 \rightarrow 0$.

For  condition (\ref{si7}), we have, by symmetry, $$\sum_{i=1}^N\sum_{j,k\in node_{G_i}}^{j\neq k}c_0^3c_j = \sum_{i=1}^N\sum_{j,k\in node_{G_i}}^{j\neq k}c_0^3c_k.$$ Thus we only need to deal with two components $\sum_{i=1}^N\sum_{j,k\in node_{G_i}}^{j\neq k}c_0^3c_i1_{\{(j,k) \text{ exists}\}}$ and \\$\sum_{i=1}^N\sum_{j,k\in node_{G_i}}^{j\neq k}c_0^3c_j$. It is not hard to see that they would go to zero under  condition $\sum_{i=1}^N|G_i|^2 = o(|G|^{1.5})$ as
\begin{align*}
    &\sum_{i=1}^N\sum_{j,k\in node_{G_i}}^{j\neq k}c_0^3c_j \leq  c_0^3\sum_{i=1}^N|G_i|\sum_{j\in node_{G_i}}c_j \leq  c_0^3\sum_{i=1}^N|G_i|\sqrt{|G_i|\sum_{j\in node_{G_i}}c_j^2}\precsim c_0^3\sum_{i=1}^N|G_i|^{1.5},\\
    &\sum_{i=1}^N\sum_{j,k\in node_{G_i}}^{j\neq k}c_0^3c_i1_{\{(j,k) \text{ exists}\}} \leq c_0^3\sum_{i=1}^Nc_i|G_i|^2 \precsim c_0^3\sum_{i=1}^N|G_i|^2.
\end{align*}

\section*{Acknowledgment}

The authors were partly supported by NSF DMS-1848579.

%%%%%%%%%%%%%%%%%%%%%%%%%%%%%%%%%%%%%%%%%%%%%%
%% Single Appendix:                         %%
%%%%%%%%%%%%%%%%%%%%%%%%%%%%%%%%%%%%%%%%%%%%%%
%\begin{appendix}
%\section*{???}%% if no title is needed, leave empty \section*{}.
%\end{appendix}
%%%%%%%%%%%%%%%%%%%%%%%%%%%%%%%%%%%%%%%%%%%%%%
%% Multiple Appendixes:                     %%
%%%%%%%%%%%%%%%%%%%%%%%%%%%%%%%%%%%%%%%%%%%%%%
%\begin{appendix}
%\section{???}
%
%\section{???}
%
%\end{appendix}

%%%%%%%%%%%%%%%%%%%%%%%%%%%%%%%%%%%%%%%%%%%%%%
%% Support information, if any,             %%
%% should be provided in the                %%
%% Acknowledgements section.                %%
%%%%%%%%%%%%%%%%%%%%%%%%%%%%%%%%%%%%%%%%%%%%%%
%\begin{acks}[Acknowledgments]
% The authors would like to thank ...
%\end{acks}
%%%%%%%%%%%%%%%%%%%%%%%%%%%%%%%%%%%%%%%%%%%%%%
%% Funding information, if any,             %%
%% should be provided in the                %%
%% funding section.                         %%
%%%%%%%%%%%%%%%%%%%%%%%%%%%%%%%%%%%%%%%%%%%%%%
%\begin{funding}
% The first author was supported by ...
%
% The second author was supported in part by ...
%\end{funding}

%%%%%%%%%%%%%%%%%%%%%%%%%%%%%%%%%%%%%%%%%%%%%%
%% Supplementary Material, including data   %%
%% sets and code, should be provided in     %%
%% {supplement} environment with title      %%
%% and short description. It cannot be      %%
%% available exclusively as external link.  %%
%% All Supplementary Material must be       %%
%% available to the reader on Project       %%
%% Euclid with the published article.       %%
%%%%%%%%%%%%%%%%%%%%%%%%%%%%%%%%%%%%%%%%%%%%%%
\begin{supplement}
\stitle{Supplementary Material to ``Limiting distributions of graph-based test statistics"}
\sdescription{It contains proofs of theorems and technical results, a comparison of the moment method and the locSCB method, and tables of the empirical size of the original, weighted, and max-type edge-count tests.}
\end{supplement}

%%%%%%%%%%%%%%%%%%%%%%%%%%%%%%%%%%%%%%%%%%%%%%%%%%%%%%%%%%%%%
%%                  The Bibliography                       %%
%%                                                         %%
%%  imsart-???.bst  will be used to                        %%
%%  create a .BBL file for submission.                     %%
%%                                                         %%
%%  Note that the displayed Bibliography will not          %%
%%  necessarily be rendered by Latex exactly as specified  %%
%%  in the online Instructions for Authors.                %%
%%                                                         %%
%%  MR numbers will be added by VTeX.                      %%
%%                                                         %%
%%  Use \cite{...} to cite references in text.             %%
%%                                                         %%
%%%%%%%%%%%%%%%%%%%%%%%%%%%%%%%%%%%%%%%%%%%%%%%%%%%%%%%%%%%%%

%% if your bibliography is in bibtex format, uncomment commands:
\bibliographystyle{imsart-nameyear} % Style BST file (imsart-number.bst or imsart-nameyear.bst)
\bibliography{reference}       % Bibliography file (usually '*.bib')

%% or include bibliography directly:
% \begin{thebibliography}{}
% \bibitem{b1}
% \end{thebibliography}

\end{document}